\documentclass[12pt]{iopart}
\usepackage{iopams}  
\usepackage{graphicx}
\newtheorem{theorem}{Theorem}[section]
\newtheorem{lemma}[theorem]{Lemma}
\newtheorem{corollary}[theorem]{Corollary}
\newtheorem{proposition}[theorem]{Proposition}

\newtheorem{assumption}[theorem]{Assumption}
\newtheorem{remark}{Remark}
\newtheorem{example}{Example}

\renewcommand{\phi}{\varphi}

\renewcommand{\ast}{\star}
\newcommand{\R}{\mathbb{R}}
\newcommand{\N}{\mathbb{N}}
\renewcommand{\bar}{\overline}
\newcommand{\calB}{\mathcal{B}}
\newcommand{\calD}{\mathcal{D}}
\newcommand{\calI}{\mathcal{I}}
\newcommand{\calR}{\mathcal{R}}
\newcommand{\calS}{\mathcal{S}}
\newcommand{\calT}{\mathcal{T}}
\newcommand{\convTX}{\stackrel{\calT_X}{\longrightarrow}}
\newcommand{\convTY}{\stackrel{\calT_Y}{\longrightarrow}}

\newcommand{\xdag}{x^\dagger}
\newcommand{\ydel}{y^\delta}
\newcommand{\xdel}{x^\delta}
\newcommand{\xMd}{\xdel_{Mo}}
\newcommand{\xId}[1]{\xdel_{#1}}
\newcommand{\xrd}{\xId{\rho}}
\newcommand{\SFyd}[1]{\calS(F(#1),\ydel)}
\newcommand{\setof}[2]{\left\{#1:#2\right\}}
\newcommand{\argmin}{\mbox{argmin}}
\newcommand{\eqref}[1]{(\ref{#1})}
\newcommand{\tfrac}[2]{{\textstyle\frac{#1}{#2}}}

\begin{document}
\title[On Ivanov and Morozov regularization]{On convergence and convergence rates for Ivanov and Morozov regularization and application to some parameter identification problems in elliptic PDEs}
\author{Barbara Kaltenbacher${^1}$ and Andrej Klassen${^2}$}
\address{${^1}$Alpen-Adria-Universit\"at Klagenfurt, Universit\"atsstra\ss e 65-67, 9020 Klagenfurt, Austria\\
${^2}$University of Duisburg-Essen, Thea-Leymann-Stra\ss e 9, 45127 Essen}
\ead{${^1}$barbara.kaltenbacher@aau.at, ${^2}$andrej.klassen@uni-due.de}
\begin{abstract}
In this paper we provide a convergence analysis of some variational methods alternative to the classical Tikhonov regularization, namely Ivanov regularization (also called method of quasi solutions) with some versions of the discrepancy principle for choosing the regularization parameter,  and Morozov regularization (also called method of the residuals). After motivating nonequivalence with Tikhonov regularization by means of an example, we prove well-definedness of the Ivanov and the Morozov method, convergence in the sense of regularization, as well as convergence rates under variational source conditions. Finally, we apply these results to some linear and nonlinear parameter identification problems in elliptic boundary value problems.
\end{abstract}
\maketitle
\section{Introduction}

Consider inverse problems formulated as operator equations
\begin{equation}\label{Fxy}
    F(x)=y
\end{equation}
$F:\calD(F)(\subseteq X)\to Y$, where $(X,\calT_X)$, $(Y,\calT_Y)$ are topological spaces. 
Such problems are typically ill-posed in the sense that $F$ is not continuously invertible.
We will assume that a solution $\xdag\in\calD(F)$ to \eqref{Fxy} exists. 

Since the measured data $\ydel$ that we actually have is typically contaminated with noise, whose level $\delta$ in the estimate
\begin{equation}\label{delta}
    \calS(y,\ydel)\leq\delta
\end{equation}
we assume to know, and due the above mentioned ill-posedness, the problem has to  be regularized.
For this purpose, we will use regularization and data misfit functionals $\calR:X\to \bar{R_0^+}$, $\calS:Y\times Y\to \bar{R_0^+}$ and consider the following two variational regularization methods.

\begin{itemize}
    \item Morozov regularization (method of the residual): $\xMd$ solves
        \begin{equation}\label{Morozov}
            \min_{x\in\calD(F)} \calR(x) \mbox{ s.t. }\calS(F(x),\ydel)\leq \tau \delta
        \end{equation}
        where $\tau\geq1$ is a fixed constant independent of $\delta$, 
cf., e.g., \cite{GrasmairHaltmeierScherzer11,LorenzWorliczek13,Morozov67} and the references therein.
    \item Ivanov regularization (method of quasi solutions): $\xrd$ solves
        \begin{equation}\label{Ivanov}
            \min_{x\in\calD(F)} \calS(F(x),\ydel) \mbox{ s.t. }\calR(x)\leq \rho
        \end{equation}
        where the radius $\rho$ of the admissible set
        \begin{equation}\label{Xadrho}
            \xdag\in X^{ad}(\rho):=\setof{x\in\calD(F)}{\calR(x)\leq\rho}
        \end{equation}
        plays the role of a regularization parameter and has to be chosen appropriately, 
cf. e.g., \cite{DombrovskajaIvanov65,Ivanov62,Ivanov63,IvanovVasinTanana02,LorenzWorliczek13,NeubauerRamlau14,SeidmanVogel89} and the references therein. 
\end{itemize}
In \cite{NeubauerRamlau14,SeidmanVogel89} well-definedness and convergence to $\xdag$ as $\delta\to0$ of $\xrd$ has been shown with the choice 
        \begin{equation}\label{rhoI}
            \rho=\rho_*^I=\rho^\dagger=\calR(\xdag)\,.
        \end{equation}
        This requires knowledge of the regularization functional value of $\xdag$, which we aim at avoiding here by choosing $\rho$ based on the discrepancy type rule
        \begin{equation}\label{rhoII}
            \eqalign{\rho=\rho_*^{II}\in\argmin\{\rho\geq0\, : &\mbox{ a minimizer }\xrd\mbox{ of \eqref{Ivanov} exists}\\&\mbox{ and } \SFyd{\xrd}\leq\tau\delta\}}
        \end{equation}
        or the relaxed, compuationally easier to fulfill version 
        \begin{equation}\label{rhoIII}
            \eqalign{\rho=\rho_*^{III}\mbox{ such that } &
                \mbox{ a minimizer }\xId{\rho_*}\mbox{ of \eqref{Ivanov} with $\rho=\rho_*$ exists and } \\
                &\delta<\SFyd{\xrd}\leq\tau\delta\,,
            }
        \end{equation}
        where like in \eqref{Morozov} $\tau\geq1$ is a fixed constant independent of $\delta$.

\medskip

As already pointed out, e.g.,  in \cite{NeubauerRamlau14,SeidmanVogel89}, the ideal choice~\eqref{rhoI} yields convergence and convergence rates even without knowledge of the noise level and of the possibly higher regularity (in the sense of a source condition) of $\xdag$. If $\calR(\xdag)$ is not known, then the discrepancy principle type choices~\eqref{rhoII}, \eqref{rhoIII} are reasonable alternatives, as will be shown here.

\medskip

We will here show convergence and convergence rates without exploiting any equvalence to Tikhonov regularization (for specially chosen regularization parameter) as this may fail due to nonconvexity, cf., e.g., \cite[Section 3.5]{IvanovVasinTanana02} and the following counterexample. 

\begin{figure}
    \includegraphics[width=0.48\textwidth]{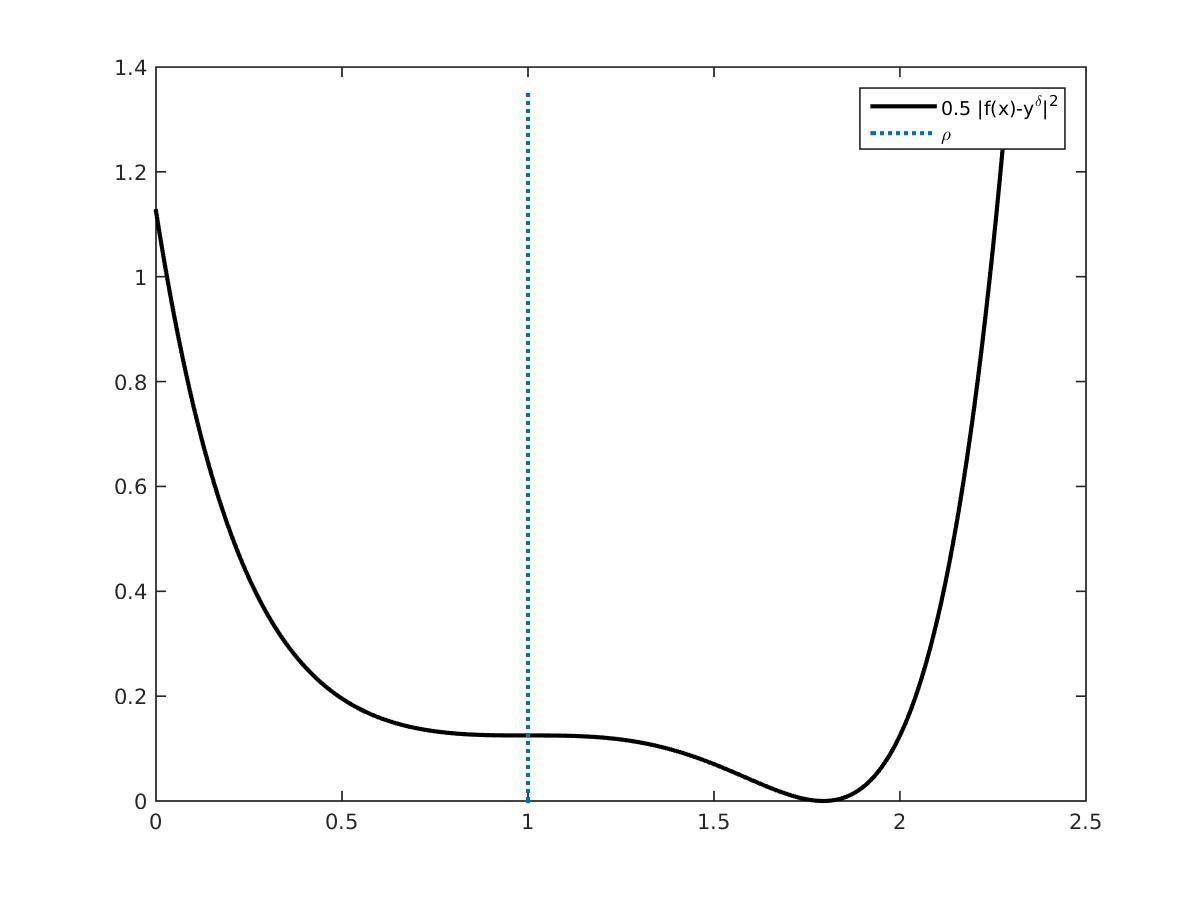}
    \hspace*{0.01\textwidth}
    \includegraphics[width=0.48\textwidth]{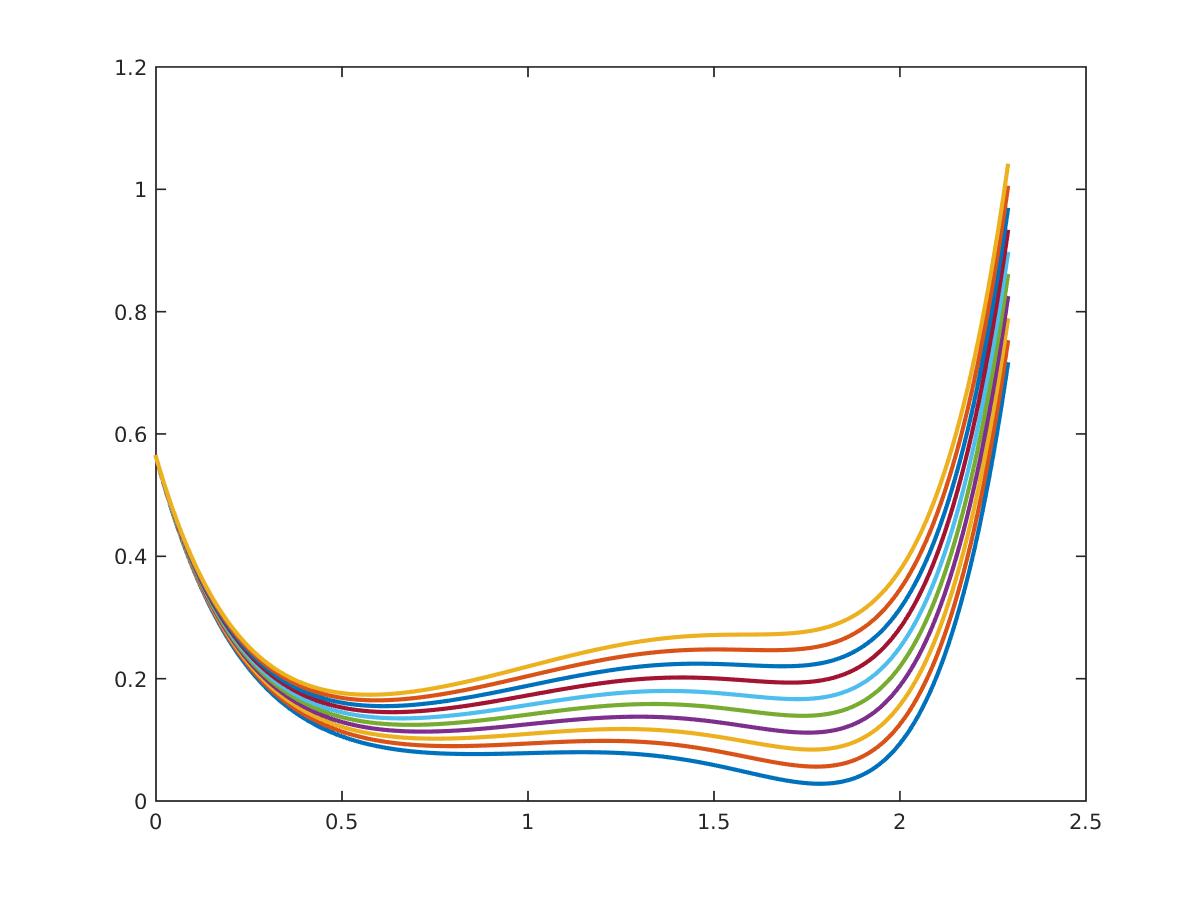}\\
    \includegraphics[width=0.48\textwidth]{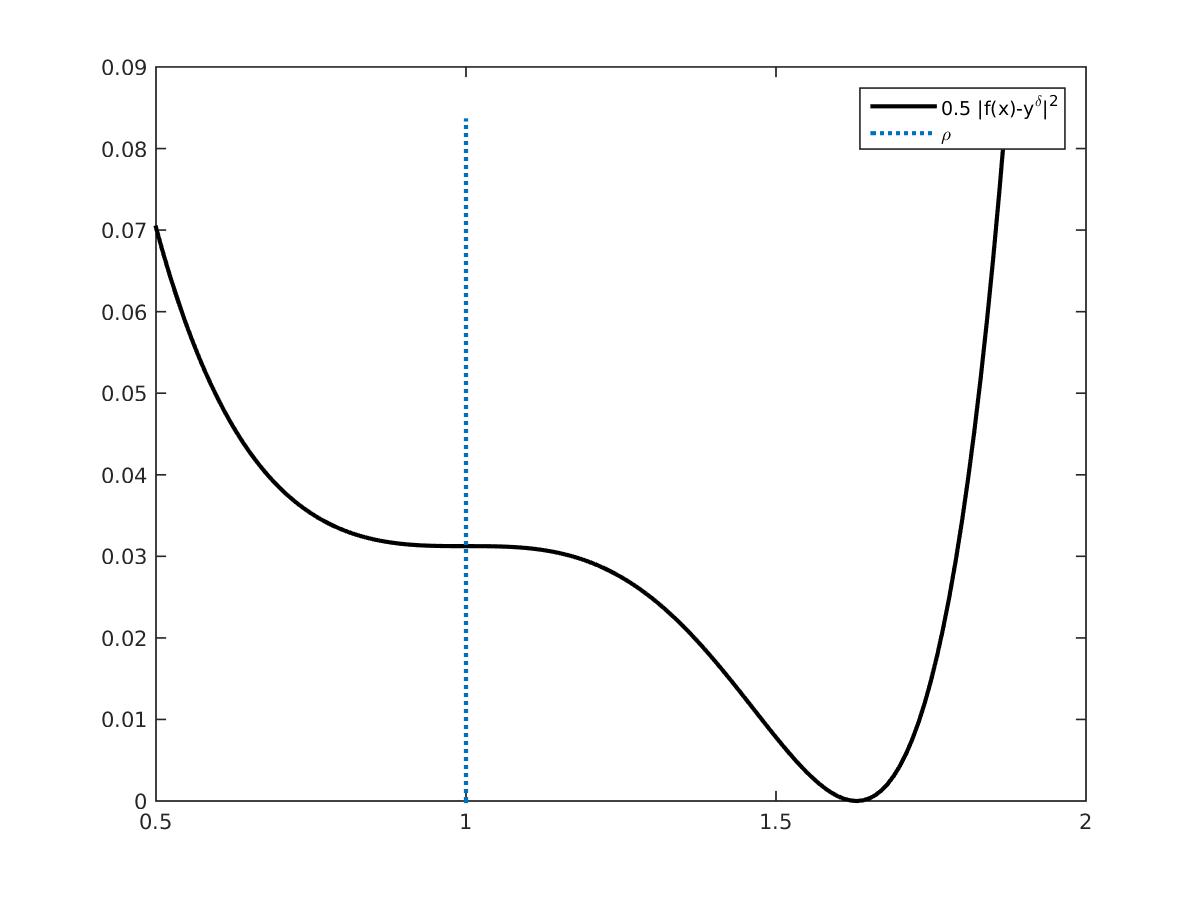}
    \hspace*{0.01\textwidth}
    \includegraphics[width=0.48\textwidth]{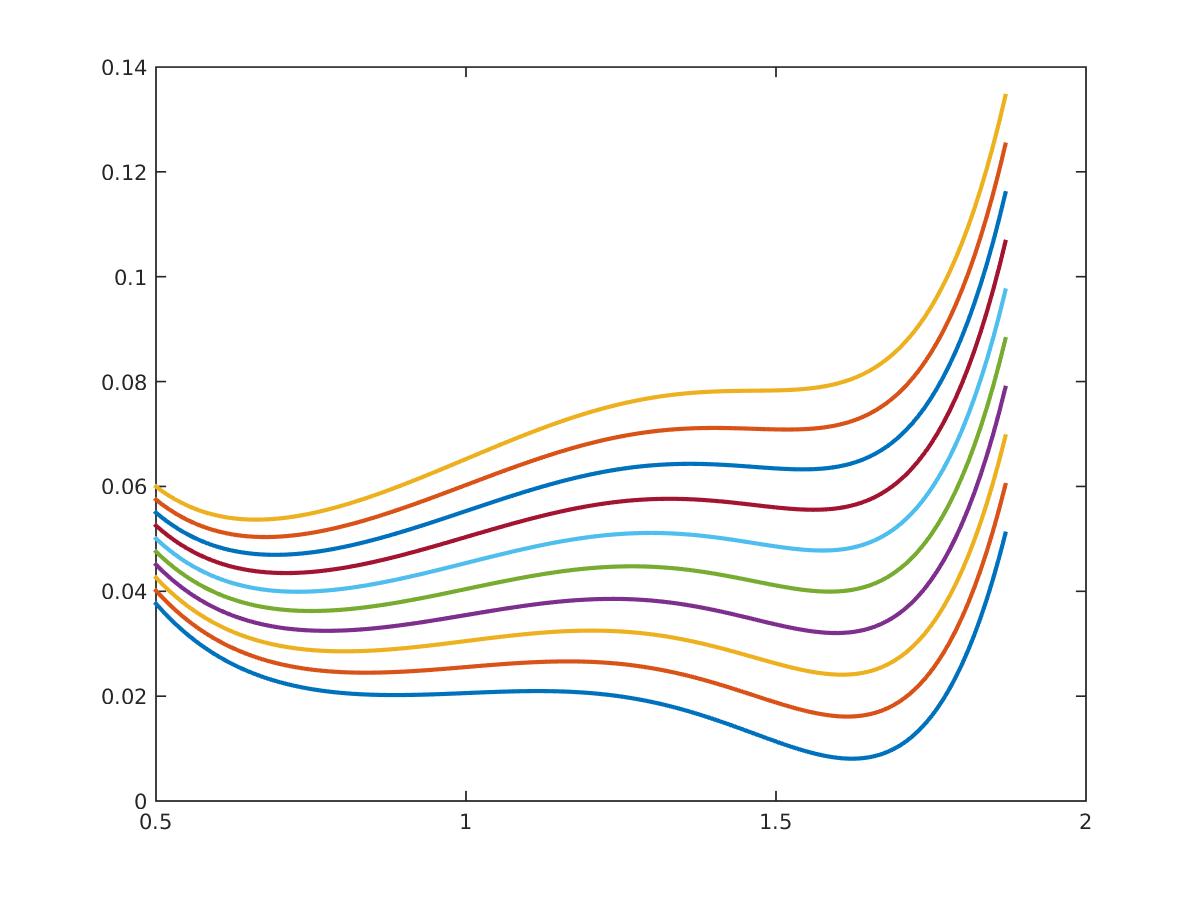}\\
    \includegraphics[width=0.48\textwidth]{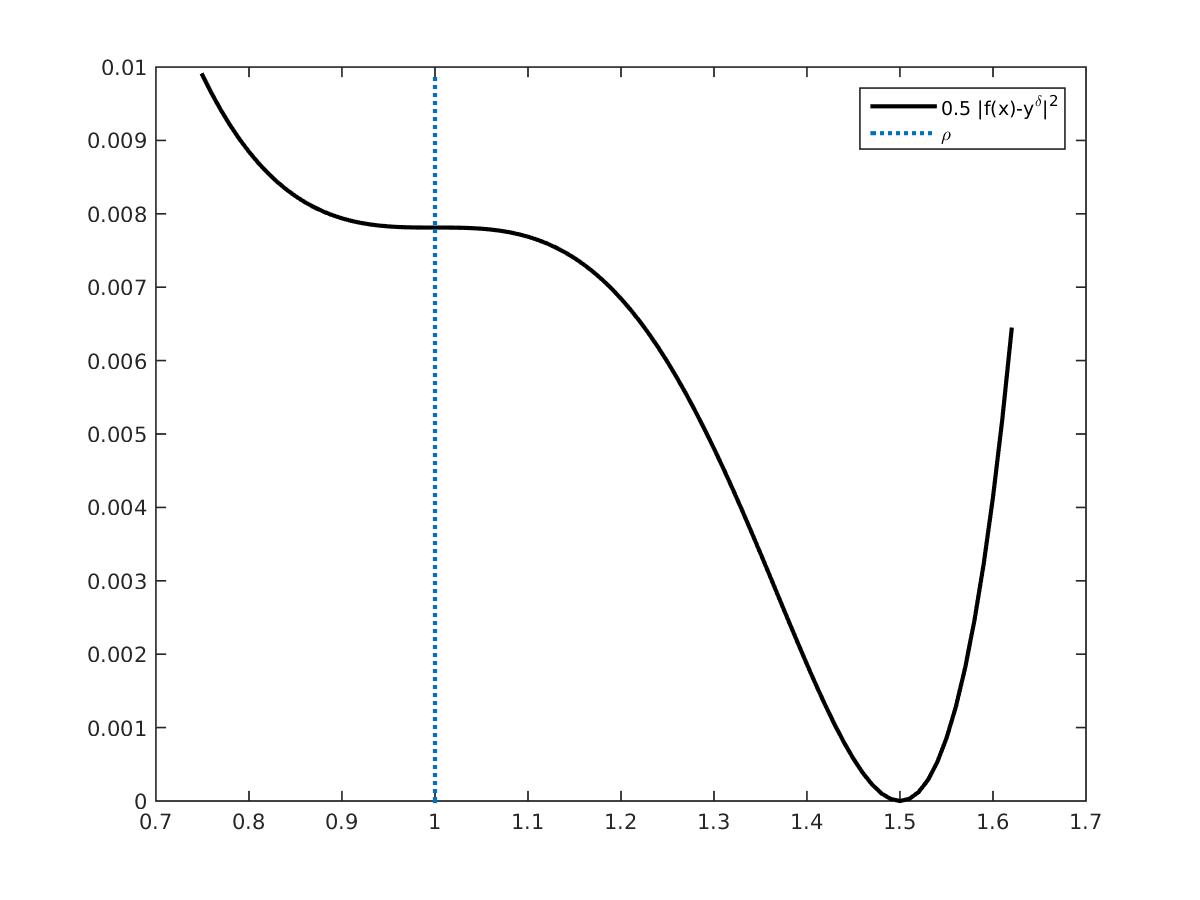}
    \hspace*{0.01\textwidth}
    \includegraphics[width=0.48\textwidth]{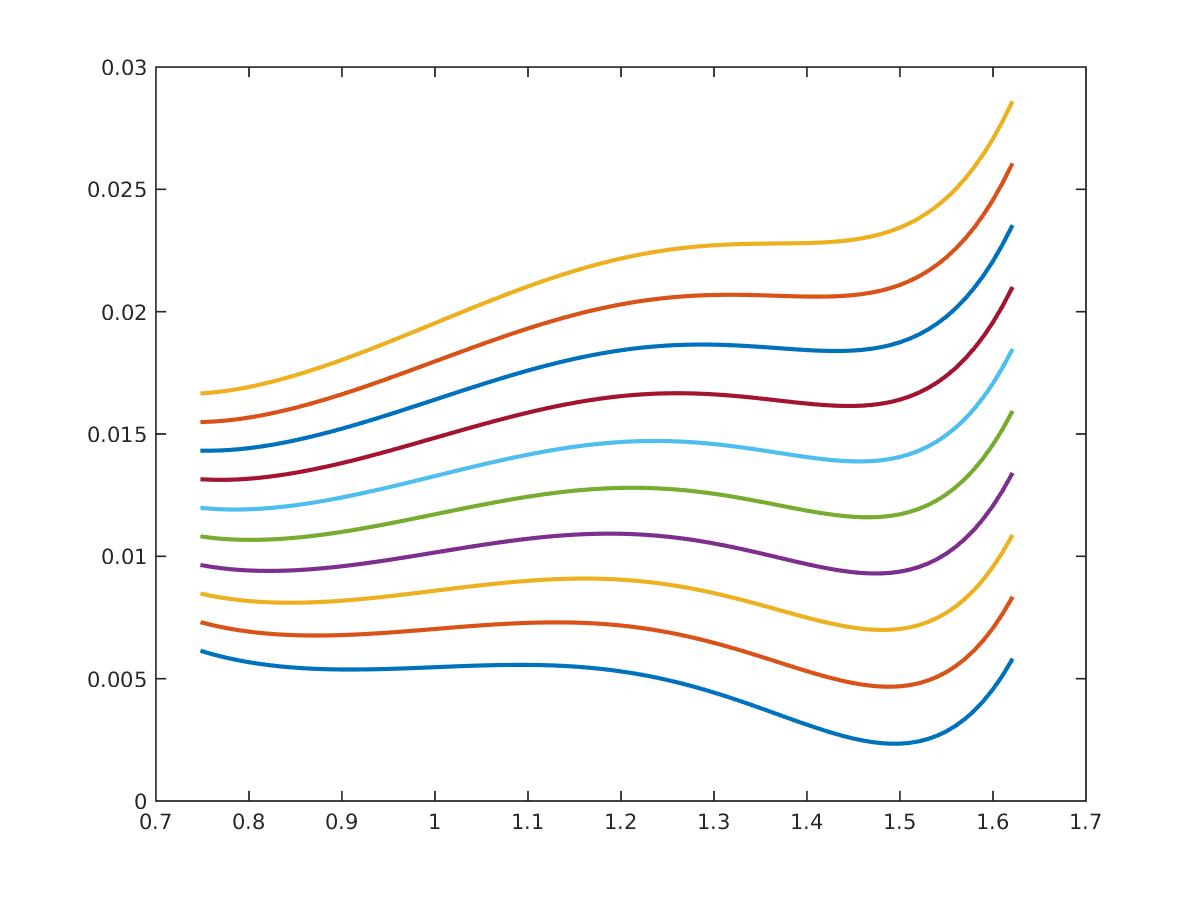}
    \caption{
        Comparison of Ivanov (left) and Tikhonov functionals (for different values of $\alpha$; right) for example \ref{counterex} with $\delta=0.5$ (top row), $\delta=0.25$ (middle row), and $\delta=0.125$ (bottom row)
        \label{fig:counterex}
    }
\end{figure}

\begin{example}\label{counterex}
    A counterexample in \cite{LorenzWorliczek13} with $F=\mbox{id}$, $x=y=\R$ and appropriately chosen $\calR$, $\calS$, shows that \eqref{Ivanov} and \eqref{Morozov} are not necessarily equivalent to Tikhonov regularization
    \begin{equation}\label{Tikhonov}
        \min_{x\in\calD(F)} \calS(F(x),\ydel) +\alpha\calR(x)\,.
    \end{equation}
    We here aim at providing a similar counterexample in the infinite dimensional Banach space setting $X=L^\infty(\Omega)$, $Y=L^2(\Omega)$ with a compact nonlinear operator $F$ and with $\calR$, $\calS$ simply chosen as powers of the respective norms.

    Consider, first of all also in the one-dimensional setting, the functionals $r(x)=|x|$, $s(x)=\frac12 |x|^2$, as well as, for fixed $\delta>0$, $x_0>(2\delta)^{1/3}>0$, $y\in\R$, $y^\delta=y+\delta$, the real function $f$ defined by 
    $f(x)=(x-x_0)^3+y$. Then $x^\dagger=x_0$ solves the exact data equation 
    $f(x)=y$ 
    as well as the Ivanov regularized problem 
    $\min_{x\in\R} |f(x)-y^\delta|$ s.t. $r(x)\leq\rho$
    with $\rho=r(x^\dagger)$ and the Morozov regularized problem 
    $\min_{x\in\R} r(x)$ s.t. $|f(x)-y^\delta|\leq\delta$.
    However, for any $\alpha>0$, a Tikonov minimizer
    $x_\alpha^\delta\in\mbox{argmin}\{ \frac12|f(x)-y^\delta|^2+\alpha r(x)\}$
    differs from $x_0$, see figure \ref{fig:counterex}. 
    
\medskip 

    This example can be lifted to an ill-posed function space setting $X=L^\infty(\Omega)$, $Y=L^2(\Omega)$, $\calS(y_1,y_2)=\frac12\|y_1-y_2\|_{L^2(\Omega)}^2$, $\calR(x)=\frac12\|x\|_{L^\infty(\Omega)}^2$, $\rho=\frac12 x_0^2$, by defining
\[
F(x)(t)=\int_\Omega\Phi(t-s)f(x(s))\, ds\,, \quad t\in\Omega
\] 
with some nonnegative normalized kernel function $\Phi:\Omega-\Omega=\setof{t-s}{s,t\in\Omega}\to[0,\infty)$ with $\int_\Omega \Phi(t-s)\, ds=1$ for all $t\in\Omega$. If $\Phi$ is weakly singular then $F$ is compact.
If $\Phi$ is the Green's function of some differential operator $D$ (equipped with boundary conditions on $\partial\Omega)$ then the operator equation $F(x)=y$ is equivalent to a possibly nonlinear inverse source problem for a PDE, namely to $Dy=f(x)$.
  
    The first order necessary optimality conditions for Ivanov and Tikhonov regularization can, analogously to Proposition 2.2 in \cite{MMC15}, and using the fact that the indicator function $\delta_{B_{x_0}^{L^\infty}(0)}$ is the Fenchel conjugate of $x\mapsto\frac{1}{x_0}\|x\|_{L^1}$, be derived as
    \begin{equation}\label{p}
        p(s)=-f'(x(s))\int_\Omega \Phi(t-s)\left(\int_\Omega \Phi(t-\tau)f(x(\tau))\, d\tau-y^\delta(t)\right)\, dt 
    \end{equation} 
    and 
    \begin{equation}\label{gradeq}
        \left\{\begin{array}{l}
            x(s)\in x_0 \mbox{sign}(p(s))\mbox{ for Ivanov regularization}\\
            x(s)\in\frac{\|p\|_{L^1}}{\alpha} \mbox{sign}(p(s))\mbox{ for Tikhonov regularization}
        \end{array}\right. \mbox{ for all } s\in\Omega
    \end{equation} 
    where 
    \[
        \mbox{sign}(a)=
        \left\{\begin{array}{ll} \{-1\} &\mbox{ if } a<0\\ \{1\} &\mbox{ if } a>0\\{} [-1,1] &\mbox{ if } a=0\end{array}\right.\,.
    \]
    Again we use constant data $\ydel(t)\equiv y+\delta$ so that by the normalization of $\Phi$ the expression for $p$ above can be rewritten as
    \[
        p(s)=-3(x(s)-x_0)^2\int_\Omega \Phi(t-s)\int_\Omega \Phi(t-\tau)((x(\tau)-x_0)^3-\delta)\, d\tau\, dt
    \]
Indeed, the only feasible element $x$ satisfying the optimality conditions for Ivanov regularization is the constant function with value $x_0$. This can be seen as follows. First of all, for any $x\in B_{x_0}^{L^\infty}(0)$, the expression $((x(\tau)-x_0)^3-\delta)$ is lower or equal $-\delta$, hence, again using the normalization of $\Phi$ we obtain $\int_\Omega \Phi(t-s)\int_\Omega \Phi(t-\tau)((x(\tau)-x_0)^3-\delta)\, d\tau\, dt\leq-\delta$, hence $p(s)\geq 3(x(s)-x_0)^2\delta\geq 0$ and $p(s)=0$ iff $x(s)=x_0$. Consider now the case $p(s)>0$, hence, by the optimality condition \eqref{gradeq}, $x(s)=x_0 \mbox{sign}(p(s))=x_0$, a contradiction.

    The functional $\Psi:x\mapsto \calS(F(x),\ydel)$ is convex on $B_{x_0}^{L^\infty}(0)$ by 
    \[
        \hspace*{-2cm}\eqalign{
            \Psi'(x)h=&\int_\Omega 3(x(s)-x_0)^2 h(s)\int_\Omega \Phi(t-s)\int_\Omega \Phi(t-\tau)((x(\tau)-x_0)^3-\delta)\, d\tau\,ds\,dt\\
            \Psi''(x)(h,h)=&\int_\Omega 6(x(s)-x_0)h^2(s)\int_\Omega \Phi(t-s)\int_\Omega \Phi(t-\tau)((x(\tau)-x_0)^3-\delta)\, d\tau\, ds\,dt\\
            &+\int_\Omega \left(\int_\Omega \Phi(t-\tau) 3 (x(\tau)-x_0)^2 h(\tau)\, d\tau\right)^2\,dt
            \quad  \geq 0
        }
    \]
    since the factors $6(x(s)-x_0)$ and $((x(\tau)-x_0)^3-\delta)$ are nonpositive for any $x\in B_{x_0}^{L^\infty}(0)$ and all the other factors are nonnegative.
    Thus we conclude from the fact that $x^\dagger:t\mapsto x_0$ obviously yields $p\equiv0$ and hence satisfies the optimality conditions $x^\dagger(t)=x_0\in[-x_0,x_0]=x_0\mbox{sign}(0)$ for Ivanov regularization, that $x^\dagger$ solves \eqref{Ivanov}. However, $x^\dagger$ does not safisfy the optimality conditions for \eqref{Tikhonov}, since $x^\dagger(t)=x_0\not\in\{0\}=\frac{\|0\|_{L^1}}{\alpha}\mbox{sign}(0)$,
    and therefore cannot be a Tikhonov minimizer.
\end{example}

Ivanov regularization has been put forward and analyzed by Ivanov and coauthors
\cite{DombrovskajaIvanov65, Ivanov62, Ivanov63, IvanovVasinTanana02}
on weakly compact sets in reflexive Banach spaces for linear inverse problems.  In \cite{SeidmanVogel89}, convergence of Tikhonov, Ivanov and Morozov regularization for nonlinear problems was established in Hilbert spaces.
More recently, a comparison of these three methods in a general setting has been provided \cite{LorenzWorliczek13} and also rates for Ivanov regularization have been established in Hilbert scales \cite{NeubauerRamlau14}. 
%A convergence analysis of Ivanov regularization for linear problems in more general Banach spaces, in particular with a nonreflexive preimage space, together with a detailed analysis of the discrepancy function is given in \cite{CK17}.
Our results on well-definedness and convergence of Morozov regularization are largely (actually in a more general framework) already covered by \cite{GrasmairHaltmeierScherzer11}, which also contains a particular convergence rates case. Nevertheless we decided to provide a joint convergence analysis with Ivanov regularization, especially in the general framework of convergence and convergence rates of Theorems \ref{th:conv}, \ref{th:rates} below, which extend the results from \cite{GrasmairHaltmeierScherzer11} also for Morozov regularization.

Parameter identification in PDEs is a class of problems, where such alternative variational formulations in Banach spaces can be particularly fruitful, e.g., when exploiting knowledge about pointwise bounds of coefficients or sources for regularization purposes. We here consider some model problems of parameter identification in the elliptic PDE
\[
\nabla\cdot (a\, \nabla u)+c\,u=b\,,
\]
namely the three possible settings of identifying one of the spatially varying parameters $a$, $b$, or $c$, from additional observations of the state $u$, while the other two parameters are assumed to be known.
For these model problems, we will establish applicability of the abstract results on Ivanov and Morozov regularization from the first part of this paper, in appropriate function space settings.

\medskip

The remainder of this paper is organized as follows.
In Section \ref{sec:conv} we provide a convergence and convergence rates analysis for \eqref{Morozov} and for \eqref{Ivanov} with \eqref{rhoI}, \eqref{rhoII}, or \eqref{rhoIII}. These abstract findings are illustrated by means of the mentioned parameter identification examples in Section \ref{sec:examples}, and we make some concluding remarks in Section \ref{sec:concl}.

\section{Convergence analysis}\label{sec:conv}
Our aim is to establish well-definedness and convergence of the methods \eqref{Morozov} and \eqref{Ivanov} with the choices \eqref{rhoI}, \eqref{rhoII}, or \eqref{rhoIII}.
For this purpose, we specify some assumptions that are actually closely related to conditions imposed in previous papers on nonlinear inverse problems. 
\begin{assumption}\label{ass1}
    There exist topologies $\calT_X$, $\calT_Y$ on $X$, $Y$ such that
    %\renenumrom
    \begin{enumerate}
        \item \label{Fclosed}
            $F$ is $\calT_X$-$\calT_Y$-sequentially closed (i.e., for any sequence $(x_n)_{n\in\N}\subseteq\calD(F)$ the implication 
            $(x_n\convTX x\mbox{ and }F(x_n)\convTY y)\ \Rightarrow \ (x\in\calD(F)\mbox{ and }F(x)=y)$
            holds.
        \item \label{RSlsc}
            $\calR$, $\calS(\cdot,\ydel)$ are 
            %convex, proper, 
            lower semicontinuous with respect to $\calT_X$ and $\calT_Y$, respectively.
        \item \label{sublevelR}
            For any $C>0$, the sublevel set $M^\calR_C=\setof{x\in\calD(F)}{\calR(x)\leq C}$ is $\calT_X$-compact. 
        \item \label{sublevelS}
            For any $C>0$ and any $\ydel\in Y$, the sublevel set $M^\calS_C=\setof{y\in F(\calD(F))}{\calS(y,\ydel)\leq C}$ is $\calT_Y$-compact.
        \item \label{Rxdag}
            There exists a solution $\xdag$ of \eqref{Fxy} such that 
            $\calR(\xdag)<\infty$. 
        \item \label{snr}
            There exists $\rho_0\geq0$ such that a minimizer $\xId{\rho_0}$ of \eqref{Ivanov} with $\rho=\rho_0$ exists and $\delta<\SFyd{\xId{\rho_0}}$.

            If $\rho_0<\rho_*^{II}$ then we additionally assume
            \begin{enumerate}
                \item \label{phin}
                    For all $\rho\in (\rho_0,\rho_*^{II}]$ and all zero sequences $(d_n)_{n\in\N}$ there exists a sequence of mappings $(\phi_n)_{n\in\N}$ such that
                    \begin{equation}
                        \eqalign{
                            &\forall n\in\N\, : \ \phi_n: X^{ad}(\rho)\to X^{ad}(\rho+d_n) \mbox{ and }\\
                            &\forall x\in X^{ad}(\rho)\, : \ \SFyd{\phi_n(x)}\to\SFyd{x} \mbox{ as }n\to\infty \,.
                        }
                    \end{equation}
                \item \label{xIunique}
                    For all $\rho\in (\rho_0,\rho_*^{II}]$ at most one minimizer $\xrd$ of \eqref{Ivanov} exists.
            \end{enumerate}
    \end{enumerate}
\end{assumption}

\begin{remark}\label{rem:examples}
    Conditions (\ref{Fclosed})-(\ref{Rxdag}) guarantee existence of an $\calR$-minimizing solution, cf. \cite[Theorem 3.4]{HKPS07}, \cite[Theorem 1.9]{PoeschlDiss}. Therefore in the following we will, without loss of generality, assume that $\xdag$ is an $\calR$-minimizing solution.

    Examples of regularization and data misfit functionals satisfying conditions (\ref{RSlsc}), (\ref{sublevelR}), (\ref{sublevelS}), (\ref{Rxdag}) are (powers of) norms on Banach spaces with the weak topology if the space is reflexive or with the weak-* topology if the space is the dual of a separable Banach space.
    With such choices of $(X,\calT_X)$, $(Y,\calT_Y)$, condition (\ref{Fclosed}) holds for any bounded linear operator $F\in L(X,Y)$.
    For some examples of nonlinear forward operators $F$ satisfying (\ref{Fclosed}), we refer to Section \ref{sec:examples}.  

    The first part of condition (\ref{snr}) is, e.g., satisfied if $\calR(x)=\|x-x_0\|^p$ and $\SFyd{x_0}>\delta$, since then by setting $\rho_0=0$ we have $\xId{\rho_0}=x_0$.

    In this norm setting $\calR(x)=\|x-x_0\|^p$, mappings $\phi_n$ according to \ref{phin} can easily be found, e.g., as $\phi_n(x):=x_0+(\tfrac{\rho+d_n}{\rho})^{1/p}(x-x_0)$, so that for any $x\in X^{ad}(\rho)$ the norm 
    $\|\phi_n(x)-x\|=(1-(\tfrac{\rho+d_n}{\rho})^{1/p})\|x-x_0\|$ tends to zero as $n\to\infty$ and \ref{phin} is satisfied if $\SFyd{\cdot}$ is continuous with respect to the norm topology.

    The uniqueness condition \ref{xIunique} is, e.g., satisfied if both functionals $\calR$ and $\SFyd{\cdot}$ are convex on $\calD(F)$ and at least one of them is strictly convex.
    In particular if $\calR$ is strictly convex, a sufficiently small radius $\rho$ might even compensate for possible nonconvexity of $F$ (analogously to sufficiently large $\alpha$ in Tikhonov regularization).

    %\item \label{SFcont}
    %$\calS\circ F$ is $\calT_Y$-continuous.
    %The continuity condition (\ref{SFcont}) is, e.g., be satisfied in case of a norm $\calS(y,\ydel)=\|y-\ydel\|$ on a Banach space $Y$ is $F$ is compact (hence, maps $\calT_X$ convergent sequences into norm convergent sequences, where $\calT_X$ might be induced by the weak or the weak-* topology, respectively, see above).

\end{remark}

    In view of the fact that the range of $F$ is typically non-closed in an ill-posed setting, the following closedness result already indicates some regularizing property of the Ivanov method. 
    \begin{proposition}\label{prop:Qrho}
        Let conditions (\ref{Fclosed}), (\ref{sublevelR}) of Assumption \ref{ass1} hold. Then for any $\rho>0$ the set $Q_\rho=F(X^{ad}(\rho))=\setof{F(x)}{x\in\calD(F)\wedge\calR(x)\leq\rho}$ is $\calT_X$-closed.
    \end{proposition}
    {\em Proof.} \
        For any sequence $(y_n)_{n\in\N}$ with $y_n\convTY y$ there exists a sequence of preimages $(x_n)_{n\in\N}\subseteq\calD(F)$ such that $\calR(x_n)\leq\rho$ and $y_n=F(x_n)$. Thus by Assumption \ref{ass1} (\ref{sublevelR}), there exists a $\calT_X$ convergent subsequence $x_{n_k}\convTX\bar{x}$ with $\calR(\bar{x})\leq\rho$, and by Assumption \ref{ass1} (\ref{Fclosed}) we get $\bar{x}\in\calD(F)$ and $F(\bar{x})=y$.
    \\ $\diamondsuit$\\

We begin our analysis with first of all showing well-definedness of minimizers.
\begin{theorem}\label{th:welldef}
    Let $\ydel\in Y$, $\tau\geq1$, $\delta>0$ be fixed and let \eqref{delta}, as well as, for two functionals $\calR:X\to \bar{R_0^+}$, $\calS:Y\times Y\to \bar{R_0^+}$, conditions (\ref{Fclosed})--(\ref{Rxdag}) of Assumption \ref{ass1} hold, 
    with $\xdag$ an $\calR$-minimizing solution of \eqref{Fxy}, 
    (which exists due to \cite[Theorem 1.9]{PoeschlDiss})

    Then $\xMd$ and, for any $\rho\geq\rho_*^{II}$, $\xrd$  are well-defined. In particular,  $\xrd$ with $\rho=\rho_*^I$ according to \eqref{rhoI} or with $\rho=\rho_*^{II}$ according to \eqref{rhoII} are well defined and the relations 
    \begin{equation}\label{rhoIIrhoI}
        \rho_*^{II}\leq \rho_*^{I}=\calR(\xdag)
    \end{equation}
    and 
    \begin{equation}\label{RSbounds1}
        \eqalign{
            &\calR(\xMd)\leq\calR(\xdag)\,, \quad
            \calR(\xId{\rho_*^{I}})\leq\calR(\xdag)\,,\quad
            \calR(\xId{\rho_*^{II}})\leq\calR(\xdag)\,,\\ 
            &\SFyd{\xMd}\leq\tau\delta\,, \quad
            \SFyd{\xId{\rho_*^{I}}}\leq\delta\,, \quad
            \SFyd{\xId{\rho_*^{II}}}\leq\tau\delta
        }
    \end{equation}
    hold.
    Moreover, the monotonicity relation
    \begin{equation}\label{mon}
        \rho_1\leq \rho_2 \ \Rightarrow \SFyd{\xId{\rho_1}}\geq \SFyd{\xId{\rho_2}}
    \end{equation}
    holds for all $\rho_1,\rho_2\geq\rho_*^{II}$ and any two minimizers $\xId{\rho_i}$ of \eqref{Ivanov} with $\rho=\rho_i$, $i\in\{1,2\}$.

    If additionally Assumption \ref{ass1} (\ref{snr}) holds, then well-definedness of $\xrd$ and the monotonicity relation \eqref{mon} extend to the interval $[\min\{\rho_0,\rho_*^{II}\},\infty)$ and also $\xrd$ with $\rho=\rho_*^{III}$ according to \eqref{rhoIII} is well defined and satisfies
    \begin{equation}\label{rhoIIIrhoI}
        \rho_*^{III}\leq \rho_*^{I}=\calR(\xdag)
    \end{equation}
    \begin{equation}\label{RSbounds2}
        \calR(\xId{\rho_*^{III}})\leq\calR(\xdag)\,,\quad
        \SFyd{\xId{\rho_*^{III}}}\leq\tau\delta
    \end{equation}

\end{theorem}

{\em Proof.} \
    The key elements of the proof are (as usual in the context of variational regularization) the direct method of calculus of variations and minimality arguments.

    {\em step 1.}
    To see existence of a minimizer of \eqref{Morozov}, note that 
    \[
        X^{ad}_{Mo}:=\setof{x\in\calD(F)}{\SFyd{x}\leq\tau\delta}
    \] 
    by \eqref{delta} contains $\xdag$ and is therefore nonempty. Since $\calR$ is nonnegative and $\calR(\xdag)$ is finite, $\calI=\inf_{x\in X^{ad}_{Mo}}\calR(x)\in[0,\calR(\xdag)]$ is finite, hence a minimizing sequence $(x_n)_{n\in\N}\subseteq X^{ad}_{Mo}$ exists such that $\calR(x_n)\to \calI$. The latter and $x_n\in X^{ad}_{Mo}$ implies boundedness of the sequences $(\calR(x_n))_{n\in\N}$, $(\SFyd{x_n})_{n\in\N}$ and thus, by Assumption \ref{ass1} (\ref{sublevelR}), (\ref{sublevelS}), existence of a subsequence and of elements $\bar{x}\in X$, $\bar{y}\in Y$ such that $x_{n_k}\convTX\bar{x}$, $F(x_{n_k})\convTX\bar{y}$, that by the weak closedness of $F$ Assumption \ref{ass1} (\ref{Fclosed}) are related by $\bar{x}\in\calD(F)$, $\bar{y}=F(\bar{x})$. By the weak lower semicontinuity of $\calR$, $\calS$ we have 
    \[
        \hspace*{-2cm}\calR(\bar{x})\leq \liminf_{k\to\infty} \calR(x_{n_k})=\calI\mbox{ and }
        \SFyd{\bar{x}}=\calS(\bar{y},\ydel)\leq \liminf_{k\to\infty} \SFyd{x_{n_k}} \leq\tau\delta\,,
    \] 
    thus $\bar{x}$ is a minimizer of \eqref{Morozov}.

    {\em step 2.} 
    Similarly, for showing existence of a minimizer of \eqref{Ivanov} with $\rho=\rho_*^I=\calR(\xdag)$ we use the fact that obviously $\xdag\in X^{ad}(\rho)$ (cf. \eqref{Xadrho})
    and therefore by \eqref{delta} and nonnegativity of $\calS$, the infimum $\calI=\inf_{x\in X^{ad}(\rho)}\SFyd{x}$ is contained in $[0,\delta]$ and thus finite. Hence, the functional values $\SFyd{x_n}$ of the minimizing sequence $(x_n)_{n\in\N}$ are bounded (by $\delta$) and boundedness of the functional values $\calR(x_n)$ follows directly from $x_n\in X^{ad}(\rho)$. The rest of the proof, using closedness of $F$ and lower semicontinuity of the functionals $\calR$, $\calS(\cdot,\xdel)$, goes analogously to above.

    {\em step 3.} well-definedness of $\rho_*^{II}$ and validity of estimates \eqref{rhoIIrhoI}, \eqref{RSbounds1}, \eqref{mon}:\\
    To prove that $\rho$ according to \eqref{rhoII} is well-defined, we show that 
    \[
        R^{ad}=\setof{\rho\geq0}{\mbox{ a minimizer }\xrd\mbox{ of \eqref{Ivanov} exists and } 
            \SFyd{\xrd}\leq\tau\delta}
    \] 
    is a right unbounded interval containing its left end point (which then is the searched for minimizer). First of all, $R^{ad}$ contains $\rho_*^I=\rho^\dagger=\calR(\xdag)$, as we have shown well-definedness of $\xId{\rho^\dagger}$ above, and since by minimality of $\xId{\rho^\dagger}$ in $X^{ad}(\rho^\dagger)\ni\xdag$ we have $\SFyd{\xId{\rho^\dagger}}\leq\SFyd{x^\dagger}\leq\delta\leq\tau\delta$. Morover for any $\underline{\rho}\in R^{ad}$ the whole interval $[\underline{\rho},\infty)$ has to be contained in $R^{ad}$, as can be easily seen by replacing $\xdag$ with $\xId{\underline{\rho}}$ in step 2. of the proof, and using the fact that $X^{ad}(\rho)\supseteq X^{ad}(\underline{\rho})$ for $\rho\geq\underline{\rho}$ so that 
    \[
        \SFyd{\xId{\rho}}=\min_{x\in X^{ad}(\rho)}\SFyd{x}\leq
        \min_{x\in X^{ad}(\underline{\rho})}\SFyd{x}= \SFyd{\xId{\underline{\rho}}}\,.
    \] 
    (By the same argument, also monotonicity \eqref{mon} follows.)
    Thus $R^{ad}$ is a union of right unbounded intervals and therefore itself a right unbounded interval. It contains its left endpoint, since actually for any sequence $\rho_n$ converging (without loss of generality in a monotonically decrasing manner) to some $\rho\in R^{ad}$, the limit $\rho$ will be contained in $R^{ad}$: Namely for all $n\in\N$, we have that $\xId{\rho_n}$ is well-defined, $\SFyd{\xId{\rho_n}}\leq\tau\delta$, $\calR(\xId{\rho_n})\leq\rho_n\leq\rho_0$, so by Assumption \ref{ass1} (\ref{sublevelR}), (\ref{sublevelS}), (\ref{Fclosed}), we have existence of a subsequence $\xId{\rho_{n_k}}$ with $\tau_X$ limit $\bar{x}$ such that $F(\xId{\rho_{n_k}})$ $\tau_Y$-converges to $F(\bar{x})$ and by the lower semicontinuity Assumption \ref{ass1} (\ref{RSlsc}) satisfies 
    \[
        \calR(\bar{x})\leq\liminf_{n\to\infty}\rho_n=\rho\,, \quad
        \SFyd{\bar{x}}\leq\liminf_{n\to\infty}\SFyd{\xId{\rho_n}}\leq\tau\delta\,.
    \]
    Thus $X^{ad}(\rho)\not=\emptyset$ and $\calI=\inf_{x\in X^{ad}(\rho)}\SFyd{x}$ is contained in $[0,\tau\delta]$ and thus finite. So existence of a minimizer $\xId{\rho}$ can be shown as in step 1., and since additionally $\bar{x}\in X^{ad}(\rho)$ implies $\SFyd{\xId{\rho}}\leq\SFyd{\bar{x}}\leq\tau\delta$, we have $\rho\in R^{ad}$.
    
    The relation \eqref{rhoIIrhoI} follows by minimality of $\rho_*^{II}$ in $R^{ad}$ and the fact that $\rho_*^I$ is contained in $R^{ad}$, as we have shown in step 2. of the proof.

    {\em step 4.} well-definedness of $\rho_*^{III}$ and validity of estimates \eqref{rhoIIIrhoI}, \eqref{RSbounds2}:\\
    Analogously to step 2. --- with $\delta$ replaced by $\SFyd{\xId{\rho_0}}$ and $\xdag$ by $\xId{\rho_0}$ --- it follows that for all $\rho\geq\min\{\rho_0,\rho_*^{II}\}$ a minimizer $\xId{\rho}$ exists. Obvioulsy also the monotonicity \eqref{mon} extends to the interval $[\min\{\rho_0,\rho_*^{II}\},\infty]$.

    If $\SFyd{\xId{\rho_0}}\leq\tau\delta$ we can set $\rho_*^{III}=\rho_0$ and are done.

    It remains to consider the case $\SFyd{\xId{\rho_0}}>\tau\delta$, which by \eqref{mon} and contraposition implies $\rho_0<\rho_*^{II}$. We will apply the Intermediate Value Theorem to the mapping $\psi:[\rho_0,\rho_*^{II}]\to\R$, $\psi(\rho)=\SFyd{\xId{\rho}}$ whose values at the endpoints satisfy $\psi(\rho_0)>\tau\delta$, $\psi(\rho_*^{II})\leq\tau\delta$, so the value $\tau\delta$ will be assumed on this interval provided $\psi$ is continuous. (Note that $\psi$ is well defined even if the optimal argument $\xId{\rho}$ is nonunique since the (globally) optimal function value is unique. However, in oder to prove continuity of the value mapping $\psi$, we will need the uniqueness assumption \ref{ass1} (\ref{xIunique}).) For arbitrary $\rho\in[\rho_0,\rho_*^{II}]$ and any sequence $(\rho_n)_{n\in\N}$ converging to $\rho$, the sequences $\calR(\xId{\rho_n})\leq\rho_n$ and $\SFyd{\xId{\rho_n}}\leq \SFyd{\xId{\rho_0}}$ (by \eqref{mon}) are bounded, thus by $\calT_X$, $\calT_Y$ compactness of sublevel sets and closedness of $F$ there exists a subsequence $\xId{\rho_{n_k}}\convTX \bar{x}$ with $F(\xId{\rho_{n_k}})\convTY F(\bar{x})$, whose limit by lower semicontinuity of $\calR$, $\calS(\cdot,\ydel)$ satisfies 
    \begin{equation}\label{RSxbar}
        \calR(\bar{x})\leq\liminf_{k\to\infty}\calR(\xId{\rho_{n_k}})\leq\rho\mbox{ and } 
        \SFyd{\bar{x}}\leq\liminf_{k\to\infty}\SFyd{\xId{\rho_{n_k}}}\,.
    \end{equation}
    Now consider an arbitrary element $x\in X^{ad}(\rho)$. Using $\phi_n$ according to assumption \ref{ass1} (\ref{phin}) (where $d_n=\rho_n-\rho$) we can render $x$ admissible for the minimization problem with radius $\rho_n$, i.e., $\phi_n(x)\in X^{ad}(\rho_n)$ and thus obtain from minimality of $\xId{\rho_n}$ that $\SFyd{\phi_n(x)}\geq \SFyd{\xId{\rho_n}}$ for all $n\in\N$. Combining this with the right hand side limit in \eqref{RSxbar} and using the fact that 
    \[
        \liminf_{k\to\infty}\SFyd{\phi_{n_k}(x)}=\SFyd{x}\,,
    \] 
    we end up with $\SFyd{x}\geq \SFyd{\bar{x}}$. Since $x\in X^{ad}$ was arbitary, the assumed uniqueness of minimizers yields $\bar{x}=\xId{\rho}$.
    Therefore, analogoulsy to \eqref{RSxbar} we get for any subsequence $(\rho_{n_m})_{m\in\N}$ of $(\rho_n)_{n\in\N}$ existence of a subsequence $(\rho_{n_{m_l}})_{l\in\N}$ such that
    \[
        \hspace*{-2cm}\eqalign{
            &\SFyd{\xId{\rho}}=\SFyd{\bar{x}}\leq\liminf_{l\to\infty}\SFyd{\xId{\rho_{n_{m_l}}}}
            \leq\limsup_{l\to\infty}\SFyd{\xId{\rho_{n_{m_l}}}}\\
            &\leq\limsup_{l\to\infty}\SFyd{\phi_{n_{m_l}}(\xId{\rho})}
            =\SFyd{\xId{\rho}}\,,
        }
    \]  
    where we have used $\phi_{n_{m_l}}(\xId{\rho})\in X^{ad}(\rho)$ in the last inequality.
    By a subsequence-subsequence argument this yields $\SFyd{\xId{\rho_n}}\to\SFyd{\xId{\rho}}$ as $n\to\infty$. 

    The estimate
    \[
        \SFyd{\xId{\rho_*^{III}}}>\delta\geq\SFyd{\xdag}\geq\SFyd{\xId{\rho^\dagger}}
    \]
    and monotonicity \eqref{mon} by contraposition yields \eqref{rhoIIIrhoI}.
\\ $\diamondsuit$\\

\begin{remark}\label{rem:stability}
    Boundedness of the $\calR$ functional values \eqref{RSbounds1} or \eqref{RSbounds2}, together with the compactness Assumption \ref{ass1} (\ref{sublevelR}) also gives subsequential type stability with respect to perturbations of the data. In case of uniqueness (Assumption \ref{ass1} \ref{xIunique}), by a subsequence-subsequence argument this yields $\calT_X$-stability.
\end{remark}

Convergence and convergence rates can be obtained from the two general results Theorems \ref{th:conv}, \ref{th:rates} below, that are quite staightforward to see.

For this purpose we need some additional assumption on $\calS$, that is obviously satisfied if $\calS$ is defined by some power of a norm. 
\begin{assumption}\label{ass2}.
    \begin{enumerate}
        \item \label{Sconv}
            For any two sequences $(y_n)_{n\in\N}$, $(\tilde{y}_n)_{n\in\N}$ we have the implication
            \[
               \hspace*{-2cm}\Bigl(\calS(y_n,\tilde{y}_n)\to0\mbox{ and } \calS(\tilde{y}_n,y)\to0 \mbox{ as }n\to\infty\Bigr)\ \Rightarrow \
                \calS(y_n,y)\to0\mbox{ as }n\to\infty\,.
            \]
        \item \label{Sdef}
            $\calS(\tilde{y},y)=0$ implies $\tilde{y}=y$
    \end{enumerate}
\end{assumption}

\begin{theorem}\label{th:conv}
    Let $y\in F(\calD(F))$ and let $(\hat{x}^\delta)_{\delta>0}\subseteq\calD(F)$ be a family of regularized approximations corresponding to a family of noisy data $(\ydel)_{\delta>0}$ satisfying \eqref{delta} such that (with $y^0:=y$) for all $\delta\geq0$ and for two functionals $\calR:X\to \bar{R_0^+}$, $\calS:Y\times Y\to \bar{R_0^+}$, 
    Assumptions \ref{ass1} (\ref{Fclosed})--(\ref{sublevelS}) and \ref{ass2} are satisfied and 
    \begin{equation}\label{RC}
        \exists C>0 \ \forall \delta>0\, : \ \calR(\hat{x}^\delta)\leq C
    \end{equation}
    and
    \begin{equation}\label{SF0delta0}
        \SFyd{\hat{x}^\delta}\to0 \mbox{ as }\delta\to0
    \end{equation}
    holds.

    Then we have $\calT_X$- subsequential convergence to a solution of \eqref{Fxy} in the sense that for any zero sequence $(\delta_n)_{n\in\N}$ the sequence $\hat{x}^{\delta_n}$ has a $\calT_X$ convergent subsequence whose limit solves \eqref{Fxy}. If the solution $\xdag$ to \eqref{Fxy} is unique, then $\hat{x}^\delta\convTX\xdag$ as $\delta\to0$.
    Moreover, if \eqref{RC} holds with $C=\calR(\xdag)$, where $\xdag$ is an $\calR$-minimizing solution of \eqref{Fxy}, then the regularization terms converge
    \begin{equation}\label{convR}
        \calR(\hat{x}^\delta)\to\calR(\xdag)\mbox{ as }\delta\to0\,.
    \end{equation}
    Thus, in case $\calR$ is a norm on a space $X$ satisfying the Kadets-Klee property, and $\calT_X$ is the weak topology on that space, altogether we even have (subsequential) norm convergence.
\end{theorem}

{\em Proof.} \
    Let $(\delta_n)_{n\in\N}$ be an arbitrary sequence converging to zero. Then by \eqref{RC}, \eqref{SF0delta0} the sequences $(\calR(\hat{x}^{\delta_n}))_{n\in\N}$, $\calS(F(\hat{x}^{\delta_n},y^{\delta_n}))_{n\in\N}$ are bounded, hence Assumption \ref{ass1} (\ref{Fclosed})--(\ref{sublevelS}) yields existence of a subsequence $(\hat{x}^{\delta_{n_k}})_{k\in\N}$ that $\calT_X$-converges to some $\bar{x}$ with 
    \[
        \calR(\bar{x})\leq \liminf_{k\to\infty} \calR(\hat{x}^{\delta_{n_k}})
    \] 
    and (using Assumption \ref{ass1} (\ref{RSlsc}) with $\delta=0$)
    \[
        \calS(F(\bar{x}),y)\leq \liminf_{k\to\infty} \calS(F(\hat{x}^{\delta_{n_k}}),y)=0\,,
    \] 
    the latter following from Assumption \ref{ass2} (\ref{Sconv}) with $y_n=F(\hat{x}^{\delta_{n_k}})$, $\tilde{y}_n=y^{\delta_n}$ and \eqref{SF0delta0}, \eqref{delta}. Thus by Assumption \ref{ass2} (\ref{Sdef}), $\bar{x}$ solves \eqref{Fxy}.

    In case of uniqueness, convergence of the whole sequence follows by a subsequence-subsequence argument.

    To show \eqref{convR}, we note that \eqref{RC}, which we here have assumed to hold with $C=\calR(\xdag)$, implies $\limsup_{\delta\to0}\calR(\hat{x}^{\delta})\leq\calR(\xdag)$.
    We now assume existence of a subsequence $\delta_n\to0$ such that $\limsup_{n\to\infty}\calR(\hat{x}^{\delta_n})<\calR(\xdag)$.
    By $\calT_X$ lower semicontinuity of $\calR$ this implies $\calR(\bar{x})<\calR(\xdag)$ for the $\calT_X$ accumulation point $\bar{x}$ whose existence we have shown above. But since $\xdag$ is an $\calR$ minimizing solution, this contradicts the fact (also proven above) that $\bar{x}$ solves \eqref{Fxy}. 
\\ $\diamondsuit$\\

\begin{corollary}\label{cor:conv}
    Let $y\in F(\calD(F))$ and let $(\ydel)_{\delta>0}$ be a family of noisy data satisfying \eqref{delta} such that (with $y^0:=y$) for all $\delta\geq0$, and for two functionals $\calR:X\to \bar{R_0^+}$, $\calS:Y\times Y\to \bar{R_0^+}$, 
    Assumptions \ref{ass1} (\ref{Fclosed})--(\ref{Rxdag}) and \ref{ass2} hold
    with $\xdag$ an $\calR$-minimizing solution of \eqref{Fxy}.
    %(which exists due to \cite[Theorem 1.9]{PoeschlDiss})
    
    Then we have $\calT_X$-subsequential convergence as $\delta\to0$ to a solution of \eqref{Fxy} for $\xMd$ and for $\xrd$ with $\rho$ according to $\eqref{rhoI}$ or \eqref{rhoII}.

    If additionally Assumption \ref{ass1} (\ref{snr}) holds, then the same holds true for $\xrd$ with $\rho$ according to \eqref{rhoIII}.
\end{corollary}

To obtain convergence rates in the Bregman distance with respect to $\calR$ 
\begin{equation}\label{Bregman}
    D_\xi(\tilde{x},x)=\calR(\tilde{x})-\calR(x)-\langle \xi, \tilde{x}-x\rangle 
\end{equation}
for some $\xi$ in the subdifferential $\partial\calR(x)$, (which is nonempty, e.g., if $\calR$ is convex) we make use of a variational source condition
\begin{equation}\label{vsc}
    \eqalign{\exists \beta\in[0,1)\  \forall \tilde{x}\in\setof{x\in\calD(F)}{\calR(\tilde{x})\leq\calR(\xdag)}\ : \\
    \hspace*{3cm}-\langle \xi^\dagger, \tilde{x}-\xdag\rangle\leq
    \beta D_{\xi^\dagger}(\tilde{x},\xdag)+\phi(\calS(F(\tilde{x}),F(\xdag))),
	}
\end{equation}
cf. \cite[Equation (14)]{HKPS07},
for some index function $\phi:\R^+\to\R^+$, (i.e., $\phi$ monotonically increasing and $\lim_{t\to0}\phi(t)=0$), and $\xi^\dagger\in\partial\calR(\xdag)$.

Moreover, Assumption \ref{ass2} (\ref{Sconv}) has to be specified as the following generalized triangle inequality.
\begin{assumption}\label{ass3}
    %\begin{enumerate}
    %\item \label{Striangle}
    There exists a constant $C_S>0$ such that for all $y_1,y_2,y_3\in Y$
    \[
        \calS(y_1,y_3)\leq C_S(\calS(y_1,y_2)+\calS(y_3,y_2))
    \]
    %\end{enumerate}
\end{assumption}

\begin{theorem}\label{th:rates}
    Let $y\in F(\calD(F))$, and let, for two functionals $\calR:X\to \bar{R_0^+}$, $\calS:Y\times Y\to \bar{R_0^+}$, Assumption \ref{ass3} as well as $\partial\calR(\xdag)\not=\emptyset$, $\calR(\xdag)<\infty$ be satisfied, where $\xdag\in\calD(F)$ is a solution to \eqref{Fxy} 
    satisfying the variational source condition \eqref{vsc}.
    Moreover, let $(\hat{x}^\delta)_{\delta>0}\subseteq\calD(F)$ be a family of regularized approximations corresponding to a family of noisy data $(\ydel)_{\delta>0}$, such that \eqref{delta} and
    \begin{equation}\label{RxRxdag}
        \forall \delta>0\, : \ \calR(\hat{x}^\delta)\leq \calR(\xdag)
    \end{equation}
    as well as
    \begin{equation}\label{SFtaudelta}
        \exists \tau>0 \ \forall \delta>0 \, : \ \SFyd{\hat{x}^\delta}\leq\tau\delta
    \end{equation}
    holds.\\
    Then $\hat{x}^\delta$ satisfies the convergence rate
    \begin{equation}\label{rate}
        D_{\xi^\dagger}(\hat{x}^\delta,\xdag)\leq \frac{1}{1-\beta}\phi(C_S(\tau+1)\,\delta)
    \end{equation}
\end{theorem}

{\em Proof.} \
    \[
        \hspace*{-2cm}\eqalign{
            &D_{\xi^\dagger}(\hat{x}^\delta,\xdag)
            =\calR(\hat{x}^\delta)-\calR(\xdag)-\langle \xi^\dagger, \hat{x}^\delta-\xdag\rangle
            \leq \beta D_{\xi^\dagger}(\hat{x}^\delta,\xdag)+\phi(\calS(F(\hat{x}^\delta),F(\xdag)))\\
            &\leq \beta D_{\xi^\dagger}(\hat{x}^\delta,\xdag)+\phi(C_S(\calS(F(\hat{x}^\delta),\ydel)+
            \calS(y,\ydel)))
            \leq \beta D_{\xi^\dagger}(\hat{x}^\delta,\xdag)+\phi(C_S(\tau+1)\delta)
        }
    \]
 $\diamondsuit$\\

\begin{remark}\label{rem:E}
    As can be easily seen, a similar result can be obtained for more general error functionals $E:X\times X\to\bar{\R_0^+}$ under a more general variational smoothness assumption (cf., e.g., \cite{BrediesLorenz09,Flemming11,Grasmair10,HofmannMathe12})
    \begin{equation}
        \exists \beta>0 \ \forall \tilde{x}\in\calD(F)\ : \ 
        \beta E(\tilde{x},\xdag)\leq
        \calR(\tilde{x})-\calR(\xdag)+\phi(\calS(F(\tilde{x}),F(\xdag)))
    \end{equation}
    or the slightly weaker condition (since we can restrict attention to elements satisfying $\calR(\tilde{x})\leq\calR(\xdag)$ and can absorb the constant $\frac{1}{\beta}$ into the function $\phi$)
    \begin{equation}\label{vsa}
        \forall \tilde{x}\in\setof{x\in\calD(F)}{\calR(\tilde{x})\leq\calR(\xdag)}\ : \ 
        E(\tilde{x},\xdag)\leq\phi(\calS(F(\tilde{x}),F(\xdag)))
    \end{equation}
    for some index function $\phi:\R^+\to\R^+$, since this by \eqref{RxRxdag}, \eqref{SFtaudelta} yields
    \[
        E(\hat{x}^\delta,\xdag)\leq \phi(C_S(\tau+1)\delta)
    \]
    A possible advantage of \eqref{vsa} is that the subdifferential of $\calR$ does not get involved and an appropriate choice of the functional $E$ might also enable to state reasonable results in the context of convex but not strictly convex $\calR$, such as the $L^1$ or the $L^\infty$ norm.
\end{remark}

\begin{corollary}\label{cor:rates}
    Let $y\in F(\calD(F))$ and let, for two functionals $\calR:X\to \bar{R_0^+}$, $\calS:Y\times Y\to \bar{R_0^+}$, Assumptions \ref{ass1} (\ref{Fclosed})--(\ref{Rxdag}) and  \ref{ass3} as well as $\partial\calR(\xdag)\not=\emptyset$ be satisfied, where $\xdag\in\calD(F)$ is a solution to \eqref{Fxy} satisfying the variational source condition \eqref{vsc}.
    Moreover, let $(\ydel)_{\delta>0}$ be a family of noisy data such that \eqref{delta} holds.

    Then the convergence rate \eqref{rate} holds for $\xMd$ and for $\xrd$ with $\rho$ according to $\eqref{rhoI}$ or \eqref{rhoII}.

    If additionally Assumption \ref{ass1} (\ref{snr}) holds then the same holds true for $\xrd$ with $\rho$ according to \eqref{rhoIII}.
\end{corollary}

\section{Examples of parameter identification problems in elliptic PDEs}\label{sec:examples}
\subsection{Identification of a source term}\label{bproblem}
We start with a linear inverse problem, namely identification of the source term $b$ in the elliptic boundary value problem 
    \begin{equation} \label{bprob} \begin{array}{rcll}
        -\Delta u &=& b&\mbox{ in }\Omega\\
        u&=&g&\mbox{ on }\partial\Omega
    \end{array}\end{equation}
    from measurements of $u$ in a smooth bounded domain $\Omega$, where $g\in H^{1/2}(\partial\Omega)$ is also given. Without loss of generality (upon subtraction of a harmonic extension of the boundary data $g$ from $b$) we can assume $g=0$. 
The forward operator 
\[
F:X=L^q(\Omega)\to Y=L^p(\Omega)\,,  \quad b\mapsto(-\Delta)^{-1}b\,,
\] 
where  $-\Delta$ is the Laplace operator with homogeneous Dirichlet boundary conditions, is well-defined and bounded, by elliptic regularity even as an operator from $L^q(\Omega)$ into $W^{2,p}(\Omega)$, provided 
\[
p,q\in[1,\infty] \,, \quad  p\leq q\,.
\]
Moreover, $F$ is linear, hence weakly closed. 
    
Now we wish to explore possible choices of distance measures $E$ satisfying \eqref{vsa} with 
    \[
        \calS(y_1,y_2)=\|y_1-y_2\|_{Y}=\|y_1-y_2\|_{L^p(\Omega)}\,, \quad
        \calR(b)= \|b\|_{L^\infty(\Omega)}\,,
    \]
under appropriate assumptions on the exact solution $b^\dagger$. 
Indeed we can estimate
\[\hspace*{-2cm}\eqalign{
&\|\tilde{b}-b^\dagger\|_{W^{-1,p}(\Omega)}
=\|(-\Delta)(\tilde{u}-u^\dagger)\|_{W^{-1,p}(\Omega)}
=\sup_{\phi\in C_c^\infty(\Omega), \|\phi\|_{W^{1,p^*}(\Omega)}\leq 1}\left|\int_\Omega \nabla (\tilde{u}-u^\dagger)\nabla \phi\, dx\right|\\
&\leq \|\nabla(\tilde{u}-u^\dagger)\|_{L^p(\Omega)}
\leq \|\tilde{u}-u^\dagger\|_{L^p(\Omega)}^{1/2} \|\tilde{u}-u^\dagger\|_{W^{2,p}(\Omega)}^{1/2}\\
&\leq \|\tilde{u}-u^\dagger\|_{L^p(\Omega)}^{1/2} 
\Bigl(\|(-\Delta)^{-1}\|_{L^q\to W^{2,p}}(\|\tilde{b}\|_{L^q(\Omega)}+\|b^\dagger\|_{L^q(\Omega)})\Bigr)^{1/2}\,,
}\]
where we have used interpolation and elliptic regularity.
The strongest case here is $p=1$, i.e. an estimate in the dual of $W^{1,1}(\Omega)$. Alternatively, by Sobolev's Lemma, for any $p>d$, $\epsilon\in(0,1-d/p]$, one gets the following result in $BV^*$
\[\hspace*{-2.5cm}\eqalign{
&\|\tilde{b}-b^\dagger\|_{BV(\Omega)^*}
=\|(-\Delta)(\tilde{u}-u^\dagger)\|_{BV(\Omega)^*}=\sup_{\phi\in C_c^\infty(\Omega), \|\phi\|_{BV(\Omega)}\leq 1}\left|\int_\Omega \nabla (\tilde{u}-u^\dagger)\nabla \phi\, dx\right|\\
&\leq \|\nabla(\tilde{u}-u^\dagger)\|_{C_0(\Omega)}
\leq C_{W^{1+d/p+\epsilon,p}\to C_0}\|\nabla(\tilde{u}-u^\dagger)\|_{W^{1+d/p+\epsilon,p}(\Omega)}\\
&\leq C_{W^{1+d/p+\epsilon,p}\to C_0}\|\tilde{u}-u^\dagger\|_{L^p(\Omega)}^{(1-d/p-\epsilon)/2} \|\tilde{u}-u^\dagger\|_{W^{2,p}(\Omega)}^{(1+d/p+\epsilon)/2}\\
&\leq C_{W^{1+d/p+\epsilon,p}\to C_0}\|\tilde{u}-u^\dagger\|_{L^p(\Omega)}^{(1-d/p-\epsilon)/2} 
\|(-\Delta)^{-1}\|_{L^p\to W^{2,p}}\Bigl(\|\tilde{b}\|_{L^p(\Omega)}+\|b^\dagger\|_{L^p(\Omega)}\Bigr)^{(1+d/p+\epsilon)/2}\,,
}\]
where the first inequality holds by definition of the BV norm.
Thus \eqref{vsa} is satisfied with 
\begin{equation}\label{EW1p}
\hspace*{-2cm}E(b_1,b_2)=\|\tilde{b}-b^\dagger\|_{W^{-1,p^*}(\Omega)} \,, \quad \varphi(t)= C \sqrt{t}
\end{equation}
or with 
\begin{equation}\label{EBV}
\hspace*{-2cm}E(b_1,b_2)=\|\tilde{b}-b^\dagger\|_{BV(\Omega)^*} \,, \quad \varphi(t)= C_\epsilon t^{(1-d/p-\epsilon)/2}
\quad \mbox{ for } p>d, \epsilon\in(0,1-d/p]
\end{equation}
with appropriate constants $C, C_\epsilon>0$.

Together with Remark \ref{rem:examples} this implies the following.
\begin{proposition}
Let $y\in W^{2,q}(\Omega)\cap H_0^1(\Omega)$, $x_0\in L^q(\Omega)$, and let $(\ydel)_{\delta\in(0,\bar{\delta}]}$ be a family of noisy data satisfying $\|y-y^\delta\|_{L^p(\Omega)}\leq\delta$ and $\|F(x_0)-y^\delta\|_{L^p(\Omega)}>\delta$.

    Then we have weak (in case of $q=\infty$, weak*) convergence in $X=L^q(\Omega)$ as $\delta\to0$ to $\xdag=-\Delta y$ for $\xMd$ and for $\xrd$ with $\rho$ according to $\eqref{rhoI}$, \eqref{rhoII}, or \eqref{rhoIII}.

Moreover, the following convergence rates hold for $\hat{x}^\delta\in\{\xMd,\xId{\rho_*^{I}},\xId{\rho_*^{II}},\xId{\rho_*^{III}}\}$
    \[ 
        \|\hat{x}^\delta-x^\dagger\|_{W^{-1,p}(\Omega)}= O(\sqrt{\delta})\,, \quad
        \|\hat{x}^\delta-x^\dagger\|_{BV(\Omega)^*}= O(\delta^{(1-\frac{d}{p}-\epsilon)/2})\,,	
    \]
the latter for any $p>d, \epsilon\in(0,1-d/p]$.
\end{proposition}

%For a computational solution strategy based on a semi-smooth Newton approach, as well as numerical results for this example we refer to \cite{CK17}.

\subsection{Identification of a potential}\label{cproblem}

    Consider identification of the spatially varying potential $c$ in the elliptic boundary value problem
    \begin{equation} \label{cprob} \begin{array}{rcll}
        -\Delta u +c u &=& f&\mbox{ in }\Omega\\
        u&=&g&\mbox{ on }\partial\Omega
    \end{array}\end{equation}
    from measurements of $u$ in $\Omega$. Here $\Omega\subseteq\R^d$ is a smooth bounded domain and $f\in H^1(\Omega)^*$, $g\in H^{1/2}(\partial\Omega)$ are given.
    We choose $X=L^\infty(\Omega)=L^1(\Omega)^*$ with $\calT_X$ the weak* topology and 
    \[
        \calD(F)=\setof{c\in L^\infty(\Omega)}{c\geq0\mbox{ a.e. }}\,,
    \]
    as well as $Y=L^p(\Omega)$ with $p\in[1,\infty]$ arbitrary and $\calT_Y$ the weak (in case $p=\infty$ weak*) topology, so that Assumption~\ref{ass1} (\ref{RSlsc}), (\ref{sublevelR}), (\ref{sublevelS}), and Assumptions~\ref{ass2}, \ref{ass3} are satisfied.

We now verify Assumption~\ref{ass1} (\ref{Fclosed}). 
    Since $\calD(F)$ is weak* closed, we have, for any sequence $(c_n)_{n\in\N}\subseteq\calD(F)$ the implication 
    $c_n\stackrel{*}{\rightharpoonup} c\ \Rightarrow \ c\in\calD(F)$. 
    It remains to show that under the additional assumption $F(c_n)\convTY u$ we also get $F(c)=u$. Denoting $u_n=F(c_n)$ and with an extension $\bar{g}\in H^1(\Omega)$ of the boundary data $g$ to $\Omega$, the weak form of \eqref{cprob} for $c=c_n$, $u=u_n$
    \begin{equation} \label{weakun} \hspace*{-2cm}\eqalign{
        &u_n-\bar{g}\in H_0^1(\Omega)\mbox{ and }\\
        &\forall \phi\in H_0^1(\Omega)\ 
        \int_\Omega \{\nabla(u_n-\bar{g})\nabla\phi + c_n (u_n-\bar{g})\phi\}\, dx 
        = \int_\Omega (f+\Delta\bar{g}-c_n\bar{g})\phi\, dx
    }\end{equation}
    with $\phi=u_n-\bar{g}$ implies
    \[\hspace*{-2cm}\eqalign{
        &\|\nabla(u_n-\bar{g})\|_{L^2(\Omega)}^2
        =\int_\Omega (f+\Delta\bar{g}-c_n\bar{g})(u_n-\bar{g})\, dx
        -\|\sqrt{c_n}(u_n-\bar{g})\|_{L^2(\Omega)}^2\\
        &\leq \|f+\Delta\bar{g}\|_{H^{-1}(\Omega)}\|\nabla(u_n-\bar{g})\|_{L^2(\Omega)}
        +\|\sqrt{c_n}\bar{g}\|_{L^2(\Omega)}\|\sqrt{c_n}(u_n-\bar{g})\|_{L^2(\Omega)}\\
        &\hspace*{9cm}-\|\sqrt{c_n}(u_n-\bar{g})\|_{L^2(\Omega)}^2\\
        &\leq \frac12 \|f+\Delta\bar{g}\|_{H^{-1}(\Omega)}^2
        +\frac12 \|\nabla(u_n-\bar{g})\|_{L^2(\Omega)}^2
        +\frac14 \|\sqrt{c_n}\bar{g}\|_{L^2(\Omega)}^2\,
    }\]
    where we have used Cauchy-Schwarz and Young's inequality as well as the norm $\|\nabla\cdot\|_{L^2(\Omega)}$ on $H_0^1(\Omega)$.
    Together with boundedness of $c_n$ in $L^\infty(\Omega)$ (following from weak* convergence and the uniform boundedness principle) we thus have uniform boundedness of $u_n$ in $H^1(\Omega)$ and thus, using compactness of the embedding $H^1(\Omega)\hookrightarrow L^2(\Omega)$, existence of a subsequence $(u_{n_k})_{k\in\N}$ and an element $\bar{u}\in H^1(\Omega)$ such that $\bar{u}-\bar{g}\in H_0^1(\Omega)$, $u_{n_k}\rightharpoonup \bar{u}$ in $H^1(\Omega)$, $u_{n_k}\to \bar{u}$ in $L^2(\Omega)$. Therewith, we get, by $u_n\convTY u$, that $\bar{u}=u$ and, using \eqref{weakun} 
    \begin{equation} \label{weaku} \hspace*{-2cm}\eqalign{
        &u-\bar{g}\in H_0^1(\Omega)\mbox{ and }\\
        &\forall \phi\in H_0^1(\Omega)\ 
        \int_\Omega \{\nabla(u-\bar{g})\nabla\phi + c (u-\bar{g})\phi\}\, dx 
        - \int_\Omega (f+\Delta\bar{g}-c\bar{g})\phi\, dx\\
        &\qquad =
        \int_\Omega \nabla(u-u_{n_k})\nabla\phi \, dx + \int_\Omega c_{n_k} (u-u_{n_k})\phi\, dx  
        + \int_\Omega(c-c_{n_k})u\phi\, dx 
    }\end{equation}
    for any $k\in\N$, where all terms on the right hand side go to zero as $k\to\infty$: The first one by $u_{n_k}\rightharpoonup u$ in $H^1(\Omega)$, the second one by boundedness of $c_{n_k}$ in $L^\infty$ and $u_{n_k}\to u$ in $L^2(\Omega)$, and the last one by  $c_n\stackrel{*}{\rightharpoonup} c$ in $L^\infty(\Omega)$ and $u\phi\in L^1(\Omega)$.
    Thus, taking the limit $k\to\infty$ in \eqref{weaku} yields $F(c)=u$.

Again we consider
     \[
        \calS(y_1,y_2)=\|y_1-y_2\|_{Y}=\|y_1-y_2\|_{L^p(\Omega)}\,, \quad
        \calR(c)= \|c\|_{L^\infty(\Omega)}\,,
    \]
and intend to find a distance measure $E$ satisfying \eqref{vsa} under appropriate assumptions on the exact solution $c^\dagger$. For this purpose we assume that the state $u^\dagger$ corresponding to the exact solution $c^\dagger$ satisfies
\begin{equation}\label{udag_cproblem1}
u^\dagger\geq \underline{u}>0\mbox{ and } u^\dagger\in W^{1,q}(\Omega)\
\mbox{ where } 
q\left\{\begin{array}{l}
=p^* \mbox{ if } p^*>d\\
\in(d,\infty] \mbox{ if } p^*=d\\
=d\mbox{ if } p^*<d
\end{array}\right.
\end{equation}
where boundedness away from zero can, e.g., be achieved by some maximum principle for the elliptic PDE \eqref{cprob} together with an assumption on nonnegativity of $f$ and positivity of $g$.

We get, similarly to subsection \ref{bproblem},
\[\eqalign{
&\|\tilde{c}-c^\dagger\|_{W^{-1,p}(\Omega)}
=\|\tfrac{1}{u^\dagger}(-\Delta+\tilde{c}\cdot\mbox{id})(\tilde{u}-u^\dagger)\|_{W^{-1,p}(\Omega)}\\
&=\sup_{\phi\in C_c^\infty(\Omega), \|\phi\|_{W^{1,p^*}(\Omega)}\leq 1}
\left|\int_\Omega \left\{\nabla (\tilde{u}-u^\dagger)\nabla \left(\tfrac{1}{u^\dagger}\phi\right)
+\frac{\tilde{c}}{u^\dagger}(\tilde{u}-u^\dagger)\phi\right\}\, dx\right|\\
&\leq C_1 \|\tilde{u}-u^\dagger\|_{W^{1,p}(\Omega)}\\
&\leq C_1 \|\tilde{u}-u^\dagger\|_{L^p(\Omega)}^{1/2} \Bigl(\|\tilde{u}\|_{W^{2,p}(\Omega)}+\|u^\dagger\|_{W^{2,p}(\Omega)}\Bigr)^{1/2}\,,
}\]
where 
\[
C_1=
\max\left\{\tfrac{1}{\underline{u}}
+\tfrac{1}{\underline{u}^2}C_{W^{1,p}\to L^r} \|\nabla u^\dagger\|_{L^q(\Omega)}\,, \
\tfrac{1}{\underline{u}} \|\tilde{c}\|_{L^\infty(\Omega)}\right\}
\]
where  
\[
r\left\{\begin{array}{l}
=\infty \mbox{ if } p^*>d\\
\in[1,\infty) \mbox{ if } p^*=d\\
=\frac{dp^*}{d-p^*}\mbox{ if } p^*<d
\end{array}\right.\,.
\]
%so that $W^{1,p^*}(\Omega)$ is continuously embedded in $L^r(\Omega)$ 
and, by elliptic regularity
\[
\|\tilde{u}\|_{W^{2,p}(\Omega)}
\leq \tilde{C} ( \|f\|_{L^p(\Omega)} + \|g\|_{W^{2-1/p,p}(\partial\Omega)})\,,
\]
and likewise for $u^\dagger$.

To estimate the $BV^*$ norm we assume, in place of \eqref{udag_cproblem1}
\begin{equation}\label{udag_cproblem2}
u^\dagger\geq \underline{u}>0\mbox{ and } u^\dagger\in W^{1,p^*}(\Omega)\cap C(\Omega)
\end{equation}
and get
\[\hspace*{-2cm}\eqalign{
&\|\tilde{c}-c^\dagger\|_{BV(\Omega)^*}
=\sup_{\phi\in C_c^\infty(\Omega), \|\phi\|_{BV(\Omega)}\leq 1}
\left|\int_\Omega \left\{\nabla (\tilde{u}-u^\dagger)\nabla \left(\tfrac{1}{u^\dagger}\phi\right)
+\frac{\tilde{c}}{u^\dagger}(\tilde{u}-u^\dagger)\phi\right\}\, dx\right|\\
&\leq \|\tfrac{1}{u^\dagger}\nabla(\tilde{u}-u^\dagger)\|_{C_0(\Omega)}
+C_{BV\to L^\infty}\left(
\frac{1}{\underline{u}^2}\|\nabla(\tilde{u}-u^\dagger)\nabla u^\dagger\|_{L^1(\Omega)}
+\frac{1}{\underline{u}}\|\tilde{c}(\tilde{u}-u^\dagger)\|_{L^1(\Omega)}\right)\\
&\leq 
\frac{1}{\underline{u}}C_{W^{1+d/p+\epsilon,p}\to C_0}\|\tilde{u}-u^\dagger\|_{L^p(\Omega)}^{(1-d/p-\epsilon)/2} 
\Bigl(\|\tilde{u}\|_{W^{2,p}(\Omega)}+\|u^\dagger\|_{W^{2,p}(\Omega)}\Bigr)^{(1+d/p+\epsilon)/2}\\
&+
C_{BV\to L^\infty}\Bigl(
\frac{1}{\underline{u}^2} \|\nabla u^\dagger\|_{L^{p^*}(\Omega)}
\|\tilde{u}-u^\dagger\|_{L^p(\Omega)}^{1/2}(\|\tilde{u}\|_{W^{2,p}(\Omega)}+\|u^\dagger\|_{W^{2,p}(\Omega)})^{1/2}\\
&\hspace*{9cm}+\frac{1}{\underline{u}}\|\tilde{c}\|_{L^\infty(\Omega)}\|\tilde{u}-u^\dagger\|_{L^1(\Omega)}\Bigr)\,,
}\]
Thus \eqref{vsa} is satisfied with \eqref{EW1p} or \eqref{EBV}.

Therewith, taking into account Remark \ref{rem:examples}, we have shown the following result.
\begin{proposition}
Let $f\in H^1(\Omega)^*$, $g\in H^{1/2}(\partial\Omega)$, $\xdag,x_0\in L^\infty(\Omega)$, let $y=u^\dagger$ satisfy \eqref{cprob} with $c=\xdag$, and let $(\ydel)_{\delta\in(0,\bar{\delta}]}$ be a family of noisy data satisfying $\|y-y^\delta\|_{L^p(\Omega)}\leq\delta$ and $\|F(x_0)-y^\delta\|_{L^p(\Omega)}>\delta$.

    Then we have weak* subsequential convergence in $X=L^\infty(\Omega)$ as $\delta\to0$ to $\xdag$ for $\xMd$ and for $\xrd$ with $\rho$ according to $\eqref{rhoI}$, \eqref{rhoII}, or \eqref{rhoIII}.

If additionally $f\in L^p(\Omega)$, $g\in W^{2-1/p,p}(\partial\Omega)$ and (a) \eqref{udag_cproblem1} or (b) \eqref{udag_cproblem2} holds, then the following convergence rates hold for $\hat{x}^\delta\in\{\xMd,\xId{\rho_*^{I}},\xId{\rho_*^{II}},\xId{\rho_*^{III}}\}$
    \[ 
        (a)\ \|\hat{x}^\delta-x^\dagger\|_{W^{-1,p}(\Omega)}= O(\sqrt{\delta})\,, \quad
        (b)\ \|\hat{x}^\delta-x^\dagger\|_{BV(\Omega)^*}= O(\delta^{(1-\frac{d}{p}-\epsilon)/2})\,,
		    \]
the latter for any $p>d, \epsilon\in(0,1-d/p]$.
\end{proposition}

\subsection{Identification of a diffusion coefficient}\label{aproblem}
    Now we consider identification of the spatially varying diffusivity $a$ in the elliptic boundary value problem
    \begin{equation} \label{aprob} \begin{array}{rcll}
        -\nabla(a\nabla u)  &=& f&\mbox{ in }\Omega\\
        u&=&g&\mbox{ on }\partial\Omega
    \end{array}\end{equation}
    from measurements of $u$ in $\Omega$. Again, $\Omega\subseteq\R^d$ is a smooth domain and $f\in H^1(\Omega)^*$, $g\in H^{1/2}(\partial\Omega)$ are given.
    Similarly to above, the domain of the forward operator 
\begin{equation}\label{Faprob}
F:\calD(F)(\subseteq X)\to Y\,,  \quad a\mapsto u \mbox{ solution to \eqref{aprob}}
\end{equation}
will be 
    \begin{equation}\label{Dfaproblem}
        \calD(F)=\setof{a\in X}{\overline{\gamma}\geq a\geq\underline{\gamma}>0\mbox{ a.e. }}\,,
    \end{equation}
    for some positive constants $\underline{\gamma}<\overline{\gamma}$, but since, as we will see below, we need a space $(X,\calT_X)$ that satisifes
    \begin{equation}\label{XcompL1}
        X\hookrightarrow L^1(\Omega)\mbox{ and $\calD(F)$ is $\calT_X$ closed} 
    \end{equation} 
    in order to obtain weak sequential closedness of $F$, the choice $X=L^\infty(\Omega)$ will not be feasible this time. However, $X=BV(\Omega)$ fulfills the requirement \eqref{XcompL1}.
    As a data space, again we use $Y=L^p(\Omega)$.  
Indeed, for 
\begin{equation}\label{paprob}
\hspace*{-2cm}(d=1\mbox{ and }p\in[1,\infty])\mbox{ or }(d=2\mbox{ and }p\in[1,\infty))\mbox{ or }d\geq3\mbox{ and }p\in[1,\frac{2d}{d-2}]),
\end{equation}
continuity of the embeddings $BV(\Omega)\to L^\infty(\Omega)$, $H^1(\Omega)\to L^p(\Omega)$ guaranteess well-definedness of $F:X\to Y$, $a\mapsto u$. Even for larger $p$, one can achieve well-definedness and a uniform $W^{1,p}(\Omega)$ bound on $F(a)$ as long as $a$ is sufficiently close to a constant.
\begin{lemma}\label{lemFa}
Let $\Omega$ be a $C^2$ domain, $p\in[2,\infty]$, $f\in (W^{1,p^*}(\Omega))^*$, $g\in W^{1-1/p,p}(\partial\Omega)$, and let $\overline{\gamma}-\underline{\gamma}$ in \eqref{Dfaproblem} be sufficiently small. Then for any $a\in \calD(F)$, the solution $u$ to \eqref{aprob} is contained in $W^{1,p}(\Omega)$ and satisfies the uniform bound
$\|u\|_{W^{1,p}(\Omega)} \leq C_{\underline{\gamma},\overline{\gamma}}$ for a constant depending only on $\underline{\gamma},\overline{\gamma}$ and the domain $\Omega$.
\end{lemma}
{\em Proof.} \
We set 
%$\bar{a}=\frac{1}{|\Omega|}\int_\Omega a dx$ 
$\bar{a}=\frac{\underline{\gamma}+\overline{\gamma}}{2}$ 
and first of all prove $W^{1,p}$ regularity of solutions to \eqref{aprob} with constant diffusion coefficient $\bar{a}$. To this end, we consider a smooth approximation $v_n+\bar{g}_n$ of $u$, satisfying \eqref{aprob} with $f,g$ replaced by $f_n\in C^\infty$, $g_n\in C^\infty$, $f_n\to f$ in $(W^{1,p^*}(\Omega))^*$, $g_n\to g$ in $W^{1-1/p,p}(\partial\Omega)$, and $\bar{g}_n$ the smooth extension of $g_n$ to the interior satisfying 
\[
\|\bar{g}_n\|_{W^{1,p}(\Omega)}\leq C_{tr} \|g_n\|_{W^{1-1/p,p}(\partial\Omega)}
\] 
according to the Trace Theorem. 
We use the Helmholtz decomposition cf. \cite[Section III.]{Galdi}  
\[
\nabla v_n|\nabla v_n|^{p-2}= \nabla \phi_n + \vec{w}_n
\] 
where
\[
\hspace*{-2cm}\phi_n\in L^1_{loc}(\Omega)\cap W^{1,p^*}(\Omega)\,, \quad \|\nabla\phi_n\|_{L^{p^*}(\Omega)}\leq C\|\nabla v_n|\nabla v_n|^{p-2}\|_{L^{p^*}(\Omega)}=C\|\nabla v_n\|_{L^p(\Omega)}^{p-1}
\]
and 
\[
\hspace*{-2cm}\vec{w}_n\in L^p(\Omega)^d\,, \quad \nabla\cdot\vec{w}_n=0 \mbox{ in }\Omega, \quad \nu\cdot\vec{w}_n=0\mbox{ on }\partial\Omega\,,
\]
and use the PDE to get the energy estimate
\[\hspace*{-2cm}\eqalign{
&\bar{a}\|\nabla v_n\|_{L^p(\Omega)}^p
= \bar{a}\int_\Omega \nabla v_n \cdot \nabla v_n|\nabla v_n|^{p-2} \, dx
= \bar{a}\int_\Omega \nabla v_n \cdot \Bigl(\nabla \phi_n + \vec{w}_n\Bigr) \, dx\\
&= \int_\Omega \bar{a}\nabla v_n \cdot \nabla \phi_n \, dx
= \int_\Omega (f_n+\bar{a}\Delta\bar{g}_n) \phi_n \, dx
\leq C_{PF}\|f_n+\bar{a}\Delta\bar{g}_n\|_{(W^{1,p^*}(\Omega))^*} \|\nabla \phi_n\|_{L^{p^*}(\Omega)}\\
&\leq C_{PF}\|f_n+\bar{a}\Delta\bar{g}_n\|_{(W^{1,p^*}(\Omega))^*} \|\nabla v_n\|_{L^p(\Omega)}^{p-1}
}\]
where we have used the Poincar\'{e}-Friedrichs inequality on $W_0^{1,p}(\Omega)$.  
Thus $\|u_n\|_{W^{1,p}(\Omega)}=\|v_n+\bar{g}_n\|_{W^{1,p}(\Omega)}$ is uniformly bounded, hence there exists a weakly (weakly${}^*$ in case of $p=\infty$) convergent subsequence, whose limit can be easily checked to coincide with $u$.\\
Thus in particular $(-\Delta_{\bar{a}})^{-1}:(W^{1,p^*}(\Omega))^*\to W_0^{1,p}(\Omega)$, mapping $f$ to a solution of \eqref{aprob} with $a=\bar{a}$, $g=0$, is bounded.\\ 
Moreover, this also implies existence of an extension $\bar{g}\in W^{1,p}(\Omega)$ of the boundary data $g$ satisfying 
    \[\eqalign{
        -\bar{a}\Delta\bar{g}  &= f\mbox{ in }\Omega\\
        \bar{g}&=g\mbox{ on }\partial\Omega
    }\]
so that \eqref{aprob} is equivalent to the fixed point equation $v=Tv$ for $v=u-\bar{g}$, where the linear operator $T:W_0^{1,p}(\Omega)\to W_0^{1,p}(\Omega)$, $Tv=(-\Delta_{\bar{a}})^{-1}[\nabla((a-\bar{a})\nabla v)$
is a contraction for $\overline{\gamma}-\underline{\gamma}$ sufficiently small:
\[\hspace*{-2.3cm}\eqalign{
&\|Tv\|_{W_0^{1,p}(\Omega)}\leq 
\|(-\Delta_{\bar{a}})^{-1}\|_{(W^{1,p^*}(\Omega))^*\to W_0^{1,p}(\Omega)}\hspace*{-0.3cm}
\sup_{\psi\in W_0^{1,p^*}(\Omega)\,, \ \|\psi\|_{W_0^{1,p^*}(\Omega)}\leq1}
\left|\int_\Omega (a-\bar{a})\nabla v \nabla \psi \, dx\right|\\
&\leq 
\|(-\Delta_{\bar{a}})^{-1}\|_{(W^{1,p^*}(\Omega))^*\to W_0^{1,p}(\Omega)}
\|a-\bar{a}\|_{L^\infty(\Omega)} \|v\|_{W_0^{1,p}(\Omega)}\\
&\leq 
\|(-\Delta_{\bar{a}})^{-1}\|_{(W^{1,p^*}(\Omega))^*\to W_0^{1,p}(\Omega)}
\tfrac{\overline{\gamma}-\underline{\gamma}}{2}  \|v\|_{W_0^{1,p}(\Omega)}\,.
}\]
$\diamondsuit$\\

We now verify weak sequential closedness of $F$.
\begin{lemma}
$F$ as defined in~\eqref{Faprob} with $X=BV(\Omega)$, $Y=L^p(\Omega)$, $f\in H^1(\Omega)^*$, $g\in H^{1/2}(\partial\Omega)$, $p$ as in \eqref{paprob}, satisfies Assumption~\ref{ass1} (\ref{Fclosed}).
\end{lemma}
{\em Proof.} \ 
We consider an arbitrary sequence $(a_n)_{n\in\N}\subseteq\calD(F)$ with $a_n\convTX a$, $u_n=F(a_n)\convTY u$, and, from \eqref{XcompL1}, immediately have $a\in\calD(F)$. So it remains to show that $F(a)=u$, which we do by using the weak form of \eqref{aprob} for $a=a_n$, $u=u_n$
    \begin{equation} \label{weakun_a} \eqalign{
        &u_n-\bar{g}\in H_0^1(\Omega)\mbox{ and }\\
        &\forall \phi\in H_0^1(\Omega)\ 
        \int_\Omega a_n\nabla u_n\nabla\phi \, dx 
        = \int_\Omega f\phi\, dx
    }\end{equation}
    where again $\bar{g}\in H^1(\Omega)$ is an extension of the boundary data $g$ to $\Omega$.
    Testing with $\phi=u_n-\bar{g}$ implies
    \[\eqalign{
        &
        %\underline{gamma}\|\nabla(u_n-\bar{g})\|_{L^2(\Omega)}^2\leq 
        \|\sqrt{a_n}\nabla(u_n-\bar{g})\|_{L^2(\Omega)}^2
        =\int_\Omega \bigl(f(u_n-\bar{g})-a_n\nabla\bar{g}\nabla(u_n-\bar{g})\bigr)\, dx\\
        &\leq \left(\frac{1}{\sqrt{\underline{\gamma}}}\|f\|_{H^{-1}(\Omega)}+\|\sqrt{a_n}\nabla\bar{g}\|_{L^2(\Omega)}\right)\|\sqrt{a_n}\nabla(u_n-\bar{g})\|_{L^2(\Omega)}
    }\]
    which by pointwise boundedness of $a_n$ from above and below implies uniform boundedness of $u_n$ in $H_{0,a}^1(\Omega)$, which we define as the closure of $C_0^\infty(\Omega)$ with respect to the norm induced by the inner product $(u,v)=\int_\Omega a \nabla u \nabla v\, dx$. Thus, by compactness of the embedding $H_{0,a}^1(\Omega)\hookrightarrow L^2(\Omega)$, we get existence of a subsequence $(u_{n_k})_{k\in\N}$ and an element $\bar{u}\in H^1(\Omega)$ such that $\bar{u}-\bar{g}\in H_0^1(\Omega)$, $u_{n_k}-\bar{g}\rightharpoonup \bar{u}-\bar{g}$ in $H_{0,a}^1(\Omega)$, $u_{n_k}\to \bar{u}$ in $L^2(\Omega)$, and, due to \eqref{XcompL1} as well as $(a_n)_{n\in\N}\subseteq \calD(F)$, also $a_{n_k}\to a$ in $L^1(\Omega)$ and $a_{n_k}\stackrel{*}\rightharpoonup a$ in $L^\infty(\Omega)$. The latter two limits imply norm convergence of $a_{n_k}$ to $a$ in $L^2(\Omega)$, since
    \begin{equation} \label{L2conv} \eqalign{
        &\|a_{n_k}-a\|_{L^2(\Omega)}=\int_\Omega(a_{n_k}-a)a_{n_k}\, dx +\int_\Omega(a_{n_k}-a)a\, dx\\ 
        &\leq \|a_{n_k}-a\|_{L^1(\Omega)}\overline{\gamma}+\langle a_{n_k}-a, a\rangle_{L^\infty,L^1}\to0\mbox{ as }k\to\infty\,.
    }\end{equation}
    Again by $u_n\convTY u$, $\bar{u}$ coincides with $u$ and, due to \eqref{weakun_a}, for any $k\in\N$ satisfies
    \begin{eqnarray} \label{weaku_a} 
        &&\hspace*{-2.3cm}u-\bar{g}\in H_0^1(\Omega)\mbox{ and }\\
        &&\hspace*{-2.3cm}\forall \phi\in C_0^\infty(\Omega)\ 
        \int_\Omega (a\nabla u\nabla\phi -f\phi)\, dx 
        =\int_\Omega a\nabla (u-u_{n_k})\nabla\phi\, dx +\int_\Omega(a-a_{n_k})\nabla u_{n_k}\nabla\phi\, dx 
    \nonumber\end{eqnarray}
    for any $k\in\N$.  Here the first term on the right hand side goes to zero as $k\to\infty$ by weak $H_{0,a}^1(\Omega)$ convergence of $u_{n_k}$ to $u$, and the second one by boundedness of $\|\nabla u_{n_k}|\|_{L^2(\Omega)}$ and \eqref{L2conv}. Due to density of $C_0^\infty$ in $H_0^1(\Omega)$, this yields $F(a)=u$.
\\ $\diamondsuit$\\

    Concerning the variational smoothness assumption \eqref{vsa} with
    \[
        \calS(y_1,y_2)=\|y_1-y_2\|_{Y}=\|y_1-y_2\|_{L^p(\Omega)}\,, \quad
        \calR(c)= \|c\|_{X}\,,
    \]
    we get the following auxiliary result
    \begin{lemma}\label{lem:vsa_a}
        Let $u^\dagger=F(a^\dagger)$ and the normed spaces  $U,V,Z$ be such that $U^*=V$ or $V^*=U$ and  
        \begin{equation}\label{denseZ_a}
            \setof{\nabla u^\dagger\cdot\nabla \phi}{\phi\in H_0^1(\Omega)}\mbox{ is dense in }Z
        \end{equation}
        \begin{equation}\label{estZ_a}
            \exists C>0\, \forall \phi\in H_0^1(\Omega)\,, a\in M_{a^\dagger}\, : 
            \ \|\nabla(a\nabla\phi)\|_V\leq C\|\nabla u^\dagger\cdot\nabla\phi\|_Z\,.
        \end{equation}
		for $M_{a^\dagger}=\setof{a\in\calD(F)}{\|a\|_{X}\leq\|a^\dagger\|_{X}}$.\\
        Then 
        \begin{equation}\label{vsa_a}
            \forall \tilde{a}\in M_{a^\dagger}\ : \ 
            \|\tilde{a}-a^\dagger\|_{Z^*}\leq C \|F(\tilde{a})-F(a^\dagger)\|_U
        \end{equation}
    \end{lemma}
    {\em Proof.} \
        With $\tilde{u}=F(\tilde{a})$ we have 
        \[\eqalign{
            &\|\tilde{a}-a^\dagger\|_{Z^*}
            =\sup_{\psi\in \setof{\nabla u^\dagger\cdot\nabla\phi}{\phi\in H_0^1(\Omega)}\,, \|\psi\|_Z=1}\left|\int_\Omega (\tilde{a}-a^\dagger)\psi\, dx\right|\\
            &=\sup_{\phi\in H_0^1(\Omega)\,, \|\nabla u^\dagger\cdot\nabla\phi\|_Z=1}
\left|\int_\Omega (\tilde{a}-a^\dagger)\nabla u^\dagger\cdot\nabla\phi\, dx \right|\\
            &=\sup_{\phi\in H_0^1(\Omega)\,, \|\nabla u^\dagger\cdot\nabla\phi\|_Z=1}
\left|\int_\Omega -\tilde{a}\nabla(\tilde{u}-u^\dagger)\nabla\phi\,dx \right|\\
            &=\sup_{\phi\in H_0^1(\Omega)\,, \|\nabla u^\dagger\cdot\nabla\phi\|_Z=1}
\left|\int_\Omega (\tilde{u}-u^\dagger) \nabla(\tilde{a}\nabla\phi)\,dx \right|
            \leq C\|\tilde{u}-u^\dagger\|_U
        }\]
    $\diamondsuit$\\
    For fixed $u^\dagger$, verification of condition \eqref{estZ_a} requires higher order regularity results for a solution $\phi$ to the transport equation $\nabla u^\dagger\cdot\nabla\phi=h$ in terms of higher order norms of $h$.
\\  
On the other hand, in case of $U=W^{s,p}(\Omega)$ one can estimate $\|F(\tilde{a})-F(a^\dagger)\|_U$ by means of an interpolation inequality. In doing so, one has to take into account that boundedness of $\tilde{a}$ in $BV$ in general does not admit a regularity result better than $F(\tilde{a})\in W^{1,p}(\Omega)$.   
\\
In one space dimension this can be done as follows. Setting $\Omega=(0,\ell)$, $s\in (\frac{1}{p},1]$, $U=W^{s,p}(0,\ell)$, $V=(W^{s,p}(0,\ell))^*$, $Z=W^{1-s,p^*}(0,\ell)=W^{1-s,p^*}_0(0,\ell)$ (where the latter identity follows from $1-s<\frac{1}{p^*}$) 
and assuming $(u^\dagger)'$ to be bounded away from zero with
\begin{equation}\label{udag_aproblem}
\tfrac{1}{(u^\dagger)'}\in C^{1-s}[0,\ell]
\end{equation}
and $f\in (W^{1,p}(0,\ell))^*$, we get, for any $\phi\in H_0^1(0,\ell)$, $a\in M_{a^\dagger}$ 
\[\eqalign{
&\|(a\phi')'\|_V
=\sup_{\psi\in C_c^\infty(0,\ell),\,\|\psi\|_{W^{s,p}(0,\ell)}\leq1} 
\left|\int_0^\ell a\phi'\psi'\, dx\right|\\
&\leq \|a\|_{BV} \sup_{\psi\in C_c^\infty(0,\ell),\,\|\psi\|_{W^{s,p}(0,\ell)}\leq1} 
\sup_{x\in(0,\ell)}\left|\int_0^x \phi'\psi'\, dx\right|\\
&\leq \|a^\dagger\|_{BV} \sup_{\tilde{\psi}\in C_c^\infty(0,\ell),\,\|\tilde{\psi}\|_{W^{s-1,p}(0,\ell)}\leq1} 
\sup_{x\in(0,\ell)}\left|\int_0^x \phi'\tilde{\psi}\, dx\right|\\
&\leq \|a^\dagger\|_{BV} \sup_{x\in(0,\ell)}\|\phi'\|_{W^{1-s,p^*}(0,x)}
= \|a^\dagger\|_{BV} \|\phi'\|_{W^{1-s,p^*}(0,\ell)}\\
&\leq \|a^\dagger\|_{BV} \|\tfrac{1}{(u^\dagger)'}\|_{C^{1-s}[0,\ell]} \|(u^\dagger)'\phi'\|_{W^{1-s,p^*}(0,\ell)}
}\]
where the first inequality follows from the definition of the BV norm.
Hence Lemma \ref{lem:vsa_a} together with interpolation implies that for all $\tilde{a}\in M_{a^\dagger}$, $\sigma=1-s\in[0,\tfrac{p-1}{p})$, $C=\|a^\dagger\|_{BV} \|\tfrac{1}{(u^\dagger)'}\|_{C^{1-s}[0,\ell]}$,
\[\eqalign{
&\|\tilde{a}-a^\dagger\|_{W^{-\sigma,p}(0,\ell)}
\leq C \|F(\tilde{a})-F(a^\dagger)\|_{W^{1-\sigma,p}(0,\ell)}\\
&\leq C \|F(\tilde{a})-F(a^\dagger)\|_{L^{p}(0,\ell)}^\sigma 
\Bigl(\|F(\tilde{a})\|_{W^{1,p}(0,\ell)}+\|F(a^\dagger)\|_{W^{1,p}(0,\ell)}\Bigr)^{1-\sigma}
}\]
where due to $f\in (W^{1,p}(0,\ell))^*$, $g\in W^{1-1/p,p}(\partial\Omega)$, and Lemma~\ref{lemFa} we have uniform boundedness of $\|F(\tilde{a})\|_{W^{1,p}(0,\ell)}$ on $M_{a^\dagger}$.

\begin{proposition}
Let $f\in H^1(\Omega)^*$, $g\in H^{1/2}(\partial\Omega)$, $\xdag,x_0\in BV(\Omega)$, let $y=u^\dagger$ satisfy \eqref{aprob} with $a=\xdag$, and let $(\ydel)_{\delta\in(0,\bar{\delta}]}$ be a family of noisy data satisfying $\|y-y^\delta\|_{L^p(\Omega)}\leq\delta$ and $\|F(x_0)-y^\delta\|_{L^p(\Omega)}>\delta$, where $p$ is as in \eqref{paprob}.

    Then we have weak* subsequential convergence in $X=BV(\Omega)$ as $\delta\to0$ to $\xdag$ for $\xMd$ and for $\xrd$ with $\rho$ according to $\eqref{rhoI}$, \eqref{rhoII}, or \eqref{rhoIII}.

If additionally $\Omega=(0,\ell)\subseteq\R^1$, $f\in (W^{1,p}(0,\ell))^*$, $g\in W^{1-1/p,p}(\partial\Omega)$, and \eqref{udag_aproblem} holds, then the following convergence rates hold for $\hat{x}^\delta\in\{\xMd,\xId{\rho_*^{I}},\xId{\rho_*^{II}},\xId{\rho_*^{III}}\}$
    \[ 
        \|\hat{x}^\delta-x^\dagger\|_{W^{-\sigma,p}(\Omega)}= O(\delta^\sigma)
		    \]
for any $\sigma\in[0,\tfrac{p-1}{p})$.
\end{proposition}

\section{Conclusions}\label{sec:concl}
In this paper we have shown convergence and convergence rates of Morozov regularization as well as of Ivanov regularization with some Morozov type regularization parameter choice strategies. Our convergence rates results are valid not only for the Bregman distance but also for other, sometimes more significant error measures, as illustrated for some examples of parameter identification problems in elliptic PDEs  here.

Regularizing by means of imposing bounds can be  a nice alternative to the Tikhonov type approach of regularizing by adding penalties, and it can lead to very efficient implementations when using appropriate optimization algorithms for the resulting constrained  minimization problems. For some numerical illustration of this fact we refer, e.g., to \cite{IvanovTRS}.
% and in \cite{minIPcomp}.

\section*{Acknowledgment}
Both authors wish to thank Christian Clason, University of Duisburg-Essen, for fruitful and inspiring discussions.\\ 
Moreover, the first author gratefully acknowledges financial support by the Austrian Science Fund FWF under the grants I2271 ``Regularization and Discretization of Inverse Problems for PDEs in Banach Spaces'' and P30054 ``Solving Inverse Problems without Forward Operators''. Part of this work has been completed during two research stays of the second author at the Alpen-Adria-Universit\" at Klagenfurt AAU, which was supported by the Karl Popper Kolleg ``Modeling-Simulation-Optimization'', funded by the AAU and by the Carin\-thian Economic Promotion Fund KWF.

\end{document}